
\input amstex
\documentstyle{amsppt}
\magnification = \magstep1
\hsize = 5.4 truein
\vsize = 8.7 truein
\baselineskip = 12pt
\NoBlackBoxes
\TagsAsMath
\define\kint{(\kappa,\ \kappa^+)}
\define\Ga{\Gamma}
\define\si{\sigma}

\define\la{\lambda}
\define\k{\kappa}

\define\a{\alpha}
\define\be{\beta}
\define\de{\delta}

\define\ga{\gamma}
\define\al{\aleph}
\define\vX{\varXi}

\define\th{\theta}
\define\N{\Cal N}
\define\M{\Cal M}
\define\om{\omega}

\define\nar{\narrower\smallskip\noindent}
\define\orh{\overset{\rightharpoonup}\to}
\define\modseq{(\Cal N_i|\ i \leq \theta)}
\define\z{\zeta}
\define\h{\chi}
\define\m{\mu}
\define\charset{\{0,\ 1,\ ?,\ !\}}
\define\lra{\longrightarrow}

\define\noi{\noindent}
%

\topmatter
\title A Combinatorial Forcing for Coding the Universe\\
by a Real when There Are No Sharps
\endtitle
\rightheadtext{Combinatorial Coding by a Real}
\author Saharon Shelah$^{1,2,4,5}$ and Lee J. Stanley$^{3,4,5}$
\endauthor
\thanks 1. Research partially supported by NSF,
the Basic Research Fund, Israel Academy of Science and MSRI \endthanks
\thanks 2.  Paper number 340 \endthanks
\thanks 3.   Research partially supported by NSF grant DMS 8806536, MSRI and
the Reidler Foundation. \endthanks
\thanks 4. We are very grateful to the referee for his rare combination of
diligence and patience, as well as for many helpful observations. \endthanks
\thanks 5. We would also like to thank many colleagues for their interest and
the administration and staff of MSRI and for their hospitality during
1989-90. \endthanks

\address Hebrew University, Rutgers University
\endaddress
\address
Lehigh University
\endaddress
\abstract
Assuming $0^\sharp$ does not exist,
we present a combinatorial approach to
Jensen's method of coding by a real.
The forcing uses combinatorial consequences of fine structure (including
the Covering Lemma, in various guises), but makes no direct appeal
to fine structure itself.
\endabstract
\endtopmatter
\subheading{\S 0. INTRODUCTION}
\bigskip
In \cite{5}, S. Friedman calls the original proof in \cite{1} of the Coding
Theorem \lq\lq one of the hardest in all of set theory.  The technical
considerations are extremely elaborate and the proof draws heavily on Jensen's
profound fine structure theory."  In addition to providing an excellent
overview, both of the general approach and of the particulars of that proof,
\cite{5} presents certain simplifications.  R. David,
\cite{2}, \cite{3}, and
in subsequent work, Friedman,
\cite{7}, \cite{8}, \cite{9}, \cite{10} and the forthcoming
\cite{12}, and David, \cite{4},
have shown how to integrate additional structure into the forcing conditions
to obtain yet stronger results.  In particular, \cite{8} has
generalized the coding method to coding over ground models where there are
measurable cardinals, while preserving the measurability of a designated
measurable cardinal (in the
extension, $V = L[\mu^*,\ R]$, where $R$ is a real
and $\mu^*$ is a normal measure on $\kappa$, extending a designated normal
measure, $\mu$, of the ground model).  Also, ignoring all references to
measures and mice, pp. 1147 - 1154 of \cite{8}
provides a the skeleton of a highly general and concise
version of coding over L (obtaining $V = L[R]$, in the extension, where $R$
is a real), with no hypotheses on the ground model (other than GCH), which
is fully developed in the forthcoming \cite{11}.

Nevertheless, the fundamental features of Jensen's
approach remain unchanged:  the obstacles (which we shall discuss shortly) in
the path of a \lq\lq naive" attempt to piece together the \lq\lq Building
Blocks" of Chapter 1 of \cite{1} are overcome by integrating fine-structural
considerations into the very definition of the forcing conditions.  As a
result, questions of uniformity, effectiveness and absoluteness of notions
involved in \lq\lq locally defined" approximations to the forcing must be faced.

Our approach will be radically different, drawing upon the great
simplification afforded by the hypothesis of the non-existence of $0^\sharp$.
One of the two main obstacles will
be overcome by a preliminary forcing.  The heart of the proof of Lemma 3,
in \cite{16},
which states that this preliminary forcing behaves as needed, involves
an appeal to the Covering Lemma.

The other main difficulty is to prove a
strategic closure property of the class of coding conditions.  This will be
done in (5.1) below; in (5.2) and (5.3) it is shown how this yields
distributivity properties of the class of coding conditions.
The material of \S 5 (and the material of
\cite{17}, upon which it draws) appeals
to strong combinatorial properties of $L$, developed in
\cite{17}.  A sketch of the results required for this paper
is presented in (1.2), (1.4), (1.5) below.
Here  again, the Covering Lemma
plays a role, this time by guaranteeing that the $L$-combinatorics give us a
handle on the situation in $V$.

The obvious downside to our approach is the
need for the hypothesis that $0^\sharp$ does not exist.
The main advantage is that the role of fine structure is \lq\lq modular":
it is crystallized in
the Covering Lemma itself, and in the L-combinatorics.  This is quite
analogous to the approach that the first author took in his proof of Strong
Covering (an early version appears as \cite{13}, while a revised version will
appear in the forthcoming \cite{14}).
Indeed, in some ways, this paper is an outgrowth of that work.
This allows for a simpler definition of
the coding conditions, involving the combinatorial apparatus, but
making no direct reference to fine-structural nor definability notions.
It is our hope that non-experts will find
this easier to use and to \lq\lq customize" for particular
applications they have in mind.

There are two other drawbacks to our approach.  The first is that
it seems to preclude obtaining sharp definability-type or minimal-degree type
of results.  In fact, this is a rather natural consequence of our
seeking a combinatorial forcing, intended for use in obtaining combinatorial
consequences.  We recognize that this is a significant departure from
the \lq\lq tradition" created by Jensen's original treatment, and that, in
the eyes of those steeped in that tradition, the entire approach may seem
somewhat unattractive.

The second is a consequence of our treatment of inaccessible
cardinals.  We treat them in a way which is much closer to our treatment of
successor cardinals than to our treatment of singular cardinals, in that, if
$\kappa$ is inaccessible, then, in any condition $p$, except for fewer than
$\kappa$ many $\alpha \in [\kappa,\ \kappa^+)$,
if $\alpha$ is mentioned in $p$, then $p$ says
{\it nothing} about a tail of the coding area for $\alpha$, whereas, if
$\kappa$ is singular, $\alpha \in [\kappa, \kappa^+)$, $\alpha$ is a multiple
of $\kappa^2$ (non-multiples of $\kappa^2$ are treated totally differently),
$\alpha$ is mentioned in $p$, then $p$ says {\it something} about a
tail of the coding area for $\alpha$.  Jensen has informed us that in early,
unpublished, versions of the coding paper (which evolved into \cite{1}), he
treated inaccessible cardinals as we do, but that he later shifted to treating
them in a way similar to his treatment of singular cardinals, in order, e.g.,
to be able to prove the preservation of small large cardinals, as in \S4.3 of
\cite{1}.  Thus, though we have not yet investigated the question,
our approach may preclude analogues of some of the results there.

We should point out that this way of dealing with inaccessibles (essentially
by requiring that the set of cardinals mentioned in a condition is an
Easton set) has two main uses.  The first has to do with
dealing with \lq\lq contamination", see below (2.3), (3.3) - (3.4), and
the portion of the \lq\lq SUMMARY AND INTRODUCTION" section, below, which
deals with these items.  The second involves the results of
\cite{17} and will be more fully discussed in (1.3).  Apparently, the
second use is really an essential feature of using ground model scales
as the main coding areas at singular cardinals
(see the discussion leading up to Lemma 5, below), while the
{\it need} for the first use is a
{\it result} of treating inaccessibles very differently than singulars.

Recently, S. Friedman has circulated a preprint
(\lq\lq A Short Proof of Jensen's Coding Theorem, Assuming Not
$0\#$) which draws, in part, on ideas of this paper and \cite{17}.

Finally, we would like to thank the referee for pointing out that the methods
of this paper are compatible with the existence of generics,
but that, unlike the more \lq\lq classical" coding methods,
this requires the techniques of \cite{6}.
\bigskip

\noindent
DISCUSSION.
\medskip

The coding theorem we prove is:
\medskip

\proclaim{Theorem 1}  If $V \models ZFC + GCH + \lq\lq 0^\sharp$ does not
exist", then there is a class forcing, $\bold P$ which preserves
ZFC, cofinalities, GCH and such that in $V^{\bold P}$, (there is a real $r$,
such that \lq\lq $V = L[r]$") holds.
\endproclaim

We should immediately point out that
the conditions, $\bold P = \bold P_{\aleph_2}$
of \S 3 add a subset $B_{\aleph_3} \subseteq \aleph_3$, rather than a real.
However, since, in the generic extension $V = L[B_{\aleph_3}]$ holds, it is an
easy matter to code $B_{\aleph_3}$ into a subset of $\aleph_2$, which, in
turn, is coded into a subset of $\aleph_1$, which, finally,
is coded into a real, using, e.g., almost disjoint coding.
It may be necessary to intersperse the forcings of \S 1.3 of \cite{1} to
reshape the intervals $(\aleph_i,\ \aleph_{i + 1})$, for $i = 0,\ 1$, but this
is not problematical, since when this is called for, the subset of
$\aleph_{i + 1}$ which we already have codes the universe.

Before embarking on the promised discussion of the obstacles to a naive
attempt to piecing together the Building Blocks of
\cite{1}, Chapter 1 (or some
variation on them), and how these obstacles are overcome in
Lemmas 3,\ 5, below,
we should note that we follow
\cite{5} for the general strategy for proving such
a coding theorem, and especially pp. 1005-1006, middle.  In particular, it
will suffice, by the arguments presented in \cite{5}, to prove the four main
properties of $\bold P$  presented there: Extendability, Distributivity,
Factoring, and Chain Condition.  The versions of these properties which we
prove reflect the differences between the detailed definition of our $\bold P$
and that considered in \cite{5}, but they are sufficiently similar that the
general arguments for their sufficiency go over to the setting of this paper.
These properties are proved in (6.1), (5.3), (4.4) and (6.2)
respectively.

The Factoring property states that for all regular $\theta >
\aleph_2$, there are $\bold P_\theta,\ \bold {\dot P}^\theta$ such that $P
\cong \bold P_\theta \ast \bold {\dot P}^\theta$.  The Distributivity property
states that $\bold P_\theta$ is $(\theta,\ \infty)$-distributive.  The Chain
Condition property states that, in $V^{\bold P_\theta},\ \bold {\dot
P}^\theta$ has the $\theta^+$-chain condition.
The proof of Distributivity in (5.3) is based on a strategic closure property
of $\bold P_\theta$ established in (5.1), together with result
of (5.2), which proves that the BAD player need not lose the game for
\lq\lq trivial" reasons.  Though the proof of (5.1) has been
rendered rather short and easy,
by the introduction of the \lq\lq very tidy" conditions,
and the preliminary results of (4.3) and (4.5),
in many ways, this result is the main lemma of the
entire paper.  It shows that, by the use of \lq\lq deactivators" and
generic scales (in addition to the ground model scales used in
setting up the main coding apparatus), we can overcome
the second of the two main obstacles to a naive attempt
to piece together the building blocks.  We turn now to a discussion of these
obstacles.

The first main obstacle simply involves the possibility of coding
$R \subseteq \kappa^+$ into a subset of
$\kappa$, when $\kappa$ is regular.  In order to use
almost disjoint set coding (or, as below, in \S\S 2,\ 3,
almost inclusion coding,
a variant used in \cite{15}, (1.3)),
we seem to need extra properties of the ground
model, or of the set R, since, in order to carry out the {\it decoding}
recursion across $[\kappa,\ \kappa^+)$ we need, e.g., an almost disjoint
sequence satisfying:
$$
\split
(\ast):\ \ \ \text{  for all } \theta\in (\kappa,\ \kappa^+),
\ (b_\alpha|\alpha \leq \theta) \in L[R \cap \theta],\\
\text{ and is \lq\lq canonically definable" there}.
\endsplit
$$
Such a $\overset{\rightharpoonup}\to{b}$ is called
{\it decodable}.
It is easy to obtain a  decodable
$\overset{\rightharpoonup}\to{b}$ if $R$ satisfies:
$$
(\ast\ast):\text{  for all  } \theta \in (\kappa,\ \kappa^+),
\ (card\ \theta)^{L[R \cap \theta]}= \kappa.
$$
If ($\ast\ast$) holds, we say that $R$ {\it promptly collapses fake cardinals}.

Of course, typically ($\ast\ast$) fails, and the
\lq\lq reshaping" conditions of 1.3 of
\cite{1}, the $F^B$ of \cite{5}, are introduced
to obtain ($\ast\ast$) in a generic
extension.  Unfortunately, the distributivity argument for the $F^B$
{\it seems} to require not merely that $H_{\gamma^+} = L_{\gamma^+}[B]$,
but that $H_{\gamma^{++}} = L_{\gamma^{++}}[B]$,
where $B \subseteq
\gamma^+$.  This will be the case if $B$ is the result of coding as far as
$\gamma^{+}$, but that is another story,
which leads to the original approach to
the Coding Theorem.

Instead, in \cite{16},
we showed, assuming GCH and that $0^{\sharp}$ does not exist:

\proclaim{Proposition 2} Let $\kappa > \aleph_1$ be a cardinal,
let $Z \subseteq
\kappa^{+\omega}$ be such that for all cardinals $\lambda$ with
$\kappa \leq \lambda \leq \kappa^{+\omega}$,
$H_\lambda = L_\lambda[Z]$.  Then, there
is a cofinality-preserving, GCH-preserving forcing,
$\bold S(\kappa)$, which adds a
$W \subseteq (\kappa,\ \kappa^{+})$\ such\ that\
$Z \in L[W,\ Z \cap \kappa]$\ and,\ for\ all\ $\kappa \leq \theta <
\kappa^+,\ (card\ \theta)^{L[W \cap \theta,\ Z \cap \kappa]}
= \kappa.$
\endproclaim

Then, starting from $\hat A \subseteq OR$ such that for all infinite cardinals
$\kappa,\ H_\kappa = L_\kappa [\hat A]$ and taking
$\bold S$ to be the product, with Easton
supports, of the $\bold S(\kappa)$ for
$\kappa = \aleph_2$ or $\kappa$ a limit cardinal, we
have:

\proclaim{Lemma 3} In $V^{\bold S}$, there is $A \subseteq OR$, such
that letting $\Lambda$ be the class of limit cardinals together with
$\aleph_2$:

$$ A = (A \cap \aleph_2) \cup
\bigcup\{A \cap (\kappa,\ \kappa^+): \kappa \in \Lambda \},
$$
such that for all infinite cardinals $\kappa,\ H_\kappa =
L_{\kappa}[A]$ and
such that for $\kappa = \aleph_2$ or $\kappa$ inaccessible, for all
$\kappa \leq \theta < \kappa^{+},
\ (card\ \theta)^{L[A \cap \theta]} = \kappa$ (for
singular $\kappa$, the last property is true with $L$ in place of
$L[A \cap \theta]$, in virtue of Covering).
\endproclaim

In virtue of the preceding discussion, we clearly have:

\proclaim{Corollary 4}  In $V^{\bold S}$, letting $A$ be as in Lemma 3, for
all regular $\kappa \geq \aleph_2$, there is decodable
$\overset{\rightharpoonup}\to{b} =
(b_\alpha|\alpha\in (\kappa,\ \kappa^+))$ of cofinal
almost disjoint subsets of $\kappa$
as above.
\endproclaim

In order to discuss the difficulty in proving the strategic closure properties
of the ${ \bold P}_\theta$, we need to say a
{\it bit} about the
coding apparatus for singular cardinals.  This material is discussed
at somewhat greater length in (1.2), (2.2) - (2.4)
and formally presented in
(3.4), (3.5), so the reader who finds
the present discussion insufficiently informative is encouraged to look ahead
to these items.

If $\kappa$ is singular and $\kappa < \alpha < \kappa^+$,
$\alpha$ a multiple of $\kappa^2$, then
the {\it main coding area} for $\alpha$
will be a cofinal subset of $\kappa$ which is the range
of a function, $f^*_\alpha$.  $f^*_\alpha$
is part of a {\it scale} between
$\kappa$ and $\kappa^+$.  The domain of
$f^*_\alpha$ is a fixed club subset, $D_\kappa$, of
the cardinals below $\kappa$, and for each $\lambda \in D_\kappa$,
$f^*_\alpha(\lambda)$ is of the form $\lambda^2\tau$, where $\tau$ is even,
$0 < \tau < \lambda^+$.  If $\kappa$ is a limit of singular cardinals, the
$\lambda \in
D_\kappa$ are all singular cardinals, while if $\kappa$ is of the form
$\mu^{+\omega}$, the $\lambda$ are all of the form
$\aleph_\tau$, $\tau > 1,\ \tau$ odd,
where $\kappa = \aleph_{\tau + \omega}$.

In a condition, $p$, in which $\kappa$
is mentioned, an initial segment,
$(\kappa,\ \delta^p(\kappa))$ of ordinals from
$(\kappa,\ \kappa^+)$ will be mentioned,
and a tail of $\lambda\in D_\kappa$ will be mentioned.  We shall require
that $\delta^p(\kappa)$ is a multiple of $\kappa^2$.
If $\kappa < \alpha < \delta^p(\kappa)$, $\alpha$ a multiple of
$\kappa^2$, then for a tail of $\lambda \in D_\kappa$,
$f^*_\alpha(\lambda)$ is mentioned in
$p$ (i.e., $f^*_\alpha(\lambda) < \delta^p(\lambda))$.
It is natural to expect, and will, in fact, be true in {\it tidy }
conditions that

{\narrower\medskip\noindent
\ $(\ast)$\
if $\delta^p(\kappa) \leq \alpha < \kappa^+$, and
$\alpha$ is a multiple of $\kappa^2$, then on a tail of $\lambda \in D_\kappa$,
$f^*_\alpha(\lambda)$ is {\it not} mentioned in $p$
(i.e. $\delta^p(\lambda) \leq f^{*}_\alpha(\lambda))$.\medskip}

If $(\ast)$ failed, then it might be impossible to extend
$p$ to a condition which mentions $\alpha$
and which \lq\lq codes correctly" at
$\alpha$, since the portion of $p$ below $\kappa$ may have
already imposed an unbounded amount of information on the
main coding area for $\alpha$.
However, $(\ast)$ is quite hard to maintain, when trying to
construct an upper bound for an increasing sequence of length $\theta =
\text{cf } \kappa$ of conditions from $\bold P_\theta$.

So, rather than require the property, we drop the requirement
that $p$ has to code correctly at all $\alpha$.  Instead, we allow
certain $\alpha$ to be \lq\lq deactivated", not used for coding.
We inherit another problem though:  how to detect deactivated ordinals.
For this we are led to introduce two auxiliary coding areas.  The first
is simply the set of multiples of $\kappa$
between $\alpha$ and $\alpha + \kappa^2$.
This area is used for coding an ordinal $h^p(\a) \geq \alpha$.
The idea is that not only $\alpha$ but all the ordinals in
$[\alpha,\ h^p(\a))$ will also be deactivated.

For {\it singular}
$\k$, we have, associated with each such $\a$, a function $\si^{p,\a}$
with domain $D_\k$ and we have that $h^p(\a)$ is the least $\gamma \geq
\a$ such that $f^*_\ga \geq^* \si^{p, \a}$;
in the notation introduced at the end of this section,
$h^p(\a) = scale(\si^{p,\a})$.  The $\si^{p,\a}$
are the {\it generic} scale functions, as opposed to the ground model
scale functions, $f^*_\a$.  We thank the referee for insisting on
the point of view that what we are really doing is forcing a generic scale
since the ground model scale is not adequate for dealing with deactivation.
In fact, for singular $\k,\ h^p(\a)$ is {\it decoded as}
$scale(\si^{p,\a})$
rather than being read directly in $s_\a$.  For   $\la$ of the form
$\aleph_\tau$ where $\tau > 1$ is odd, we also have $h^p(\eta)$ for multiples,
$\eta$, of $\la^2$ which are mentioned in $p$.  Here, however, there is no
associated function and the $h^p(\a)$ are directly decoded from $s_\a$.

When $\alpha$ is a limit of multiples of $\kappa^2$, the second auxiliary
coding area will be a club subset, $C_\alpha \subseteq
\alpha$.  This will be used to help us detect deactivation.  The
$C_\alpha$ will be part of a \lq\lq square system" between
$\kappa$ and $\kappa^2$.

Returning to $(\ast)$, we have mentioned that we do require it in
{\it tidy } conditions, and we require something even stronger in
{\it very tidy } conditions.  In (4.3), we show that the latter are dense.
However, by dropping the requirement $(\ast)$, we make it easier to
construct upper bounds which might not be tidy, as we do in (4.5).

In (1.1), we define games
$G(\theta,\ \orh {\Cal N},\ p_0)$, where
$\theta > \aleph_1$ is regular,
$\orh {\Cal N}$ is a certain kind of
sequence (of length $\theta + 1$) of models,
and $p_0 \in P_\theta$.
The two players, GOOD and BAD,
alternately pick conditions, $p_i \in P_\theta$.  GOOD plays at non-zero
even stages, and BAD plays at odd stages.  BAD also picks a subsequence
of $\orh {\Cal N}$ by choosing an increasing sequence $(\alpha(i)|i < \theta)$
from $\theta$; of course, $\alpha(i)$ is chosen at stage $2i + 1$.
We require that the $p_i$ are increasing,
and that $p_{2i},
\ p_{2i + 1} \in |\Cal N_{\alpha(i)}|$.  In the cases of interest,
$\orh {\Cal N}$ will satisfy a more technical
condition, introduced in (1.3),
called supercoherence.
This will guarantee that at limit stages, we will have the hypotheses
of (4.5).  GOOD wins if she succeeds in playing $p_\theta$.
BAD wins if at some even stage $j \leq \theta$, GOOD has no legal move.

In (5.1) we prove:

\proclaim{Lemma 5} For $\theta,\ p_0$ as above, and for
supercoherent $\orh {\Cal N}$,
GOOD has a winning strategy in $G(\theta,\ \orh {\Cal N},\ p_0)$.
\endproclaim

\noindent
In (5.3), it is argued (using results of (5.2) and
\cite{17}) that this gives that
$\bold P_\theta$ is
$(\theta,\ \infty$)-distributive.
What is at issue here is whether BAD always loses
because of his inability to play super-coherent sequences.
The results of \cite{17}, summarized in (1.4), below,
show that this is {\bf not} the case:  there are enough
supercoherent sequences.  In (1.4), this is presented as
a property of the combined squares and scales system, introduced in (1.2).
\bigskip
\noindent
SUMMARY AND ORGANIZATION.
\medskip

In \S 1, we present the coding apparatus for singular cardinals
and the related results from \cite{17}
notably the result about the existence of
supercoherent sequences.  In (1.1), we introduce the models sequences and the
games $G(\theta,\ \orh {\Cal N},\ p_0)$.  In (1.2) we
introduce the combinatorial apparatus
of squares and scales.  In (1.3), we introduce
the notion of supercoherence. In (1.4) we state the main result
of \cite{17}, presented as an additional property of the combinatorial
apparatus.  In (1.5), we state a
small combinatorial result about the system of scales which we use in
(4.3).  The result is clearly closely related to the definition of
very tidy condition.  This is also proved in \cite{17}.

\S 2 takes care of some other preliminaries.  In (2.1), we recapitulate
some of the material of Lemma 3, above, and (1.2), by giving a complete
discussion of coding areas for various kinds of ordinals.  In particular,
in (2.1.1) we cite an additional result from \cite{17} which shows that,
without loss of generality, we can assume that the system of $b_\eta$, for
$\eta$ such that $card\ eta$ is inaccessible has an additional property
called {\it tree-like}.  In (2.2),
we give a preliminary idea of the nature of conditions, by
introducing the class of \lq\lq protoconditions, $P(0)$.
In (2.3), we discuss the phenomenon of \lq\lq contamination"
at limit cardinals, and the devices for dealing with it,
namely the sets $X_\ga$ of \lq\lq candidates" for coding $\ga$ which are
not multiples of $\k$.
We also introduce the weak deactivator, ?, and the component,
$\beta^p$,
of conditions which provides bounds for contamination.  The $X_\ga$
are also useful in the context of the strong deactivator, !,  which we discuss
in (2.4), along with the generic scales.
In (2.5), we give a very brief sketch of the decoding procedure,
which we complete in (4.6).

In \S 3, we give the formal definition of the class of
coding conditions.  (3.1) recalls some notation, terminology and
conventions.  In (3.2) we define a sub-class
$\tilde P$ of the \lq\lq proto-conditions" of (2.2).  These still
incorporate none of the sophisticated
properties intended to deal with contamination and deactivators.  In (3.3),
we formally define the notions associated with contamination, and in
(3.4) we cut down $\tilde P$ still further by imposing five additional
properties.  The first four of these deal with contamination.  The last
deals with the use of the auxiliary coding areas,
$s_\alpha$, and thus foreshadows
(3.5), where we deal with
the strong deactivator, !, the generic scales,
$\si^{p,\a}$ and finally define
the class of coding conditions by imposing four additional properties related
to these.  In (3.6) we give the (very simple) definition of the partial
ordering of conditions.

In \S 4, we prove some basic Lemmas which will greatly facilitate our work
in \S\S 5 and 6.  In (4.1) we introduce the tidy and very tidy conditions.
We develop some of their properties in (4.1) and (4.2), and in (4.3)
we prove the crucial result that the very tidy conditions are dense.
In (4.4) we develop the Factoring Property.  In (4.5) we show that certain
increasing sequences have least upper bounds.  Taken together, (4.3) and
(4.5) do most of the groundwork for (5.1).
In (4.6), we provide a fully detailed
discussion of the decoding procedure, completing the sketch of (2.4).

In \S 5, we first prove Lemma 5, above, in (5.1).  In (5.2),
we show that the results of
\cite{17} really do mean that BAD need not
lose for trivial reasons, and in (5.3), we show that this yields the
$(\theta,\ \infty)$-distributivity of $\bold P_\theta$.  We close with two
remarks, in (5.4).  The first has to do with iterations of
$\bold P_\theta$.
The second concerns a variant of the games
$G(\theta,\ \orh {\Cal N},\ p_0)$, which we use in case (c) of (6.1)(7).
In (6.1), we establish the Extendability properties of $\bold P$, and in (6.2)
we establish the Chain Condition property.
\bigskip

\noindent
NOTATION AND TERMINOLOGY.
\medskip
Our notation and terminology is intended to be standard, or have a clear
meaning, e.g., $o.t.$ for order type, $card$ for cardinality.  A catalogue of
possible exceptions follows.  Also, the index of notation at the end of this
section summarize what follows but also some of the important definitions
and notation which is introduced in later sections.
When forcing, $p \leq q$ means $q$ gives more
information.  Closed unbounded sets are clubs.  The set of limit points of a
set $X$ of ordinals is denoted by $X^\prime$.  $A \Delta B$ is the symmetric
difference of $A$ and $B$, and $A\setminus B$ is the relative complement of
$B$ in $A$.  Notions like $=,\ \leq,\ \subseteq$, etc., when decorated
with a superscript *, mean \lq\lq on a tail".
For ordinals, $\alpha \leq \beta$, $[\alpha,\ \beta)$  is the
half-open interval $\{\gamma: \alpha \leq \gamma < \beta\}$.
The notation for the other three
intervals is clear.  It should be clear from context whether the open
interval or the ordered pair is meant.  For ordinals $\alpha,\
\beta$, we write $\alpha >> \beta$ to mean that $\alpha$ is MUCH greater than
$\beta$; the precise sense of how much greater we must take it to be
is supposed to be clear from context.

For infinite cardinals, $\kappa,\ H_\kappa$ is
the set of all sets hereditarily of cardinality $< \kappa$, i.e. those sets $x$
such that if $t$ is the transitive closure of $x$, then $card\ t <
\kappa$.  We regard $\omega$ as a successor cardinal, by ignoring the positive
finite cardinals.  Thus, for us, $\omega = 0^+$.
We say that a cardinal $\kappa$ is {\bf s-like }
if it singular or of the form $\aleph_\tau$ where $\tau > 1$ is odd,
and that it is {\bf i-like} if is inaccessible or
of the form $\aleph_\tau$ where
$\tau > 0$ is even.  For s-like cardinals,
$\kappa$, we define $U(\k)$ to be the
set of multiples of $\k^2$ in $(\kappa,\ \kappa^+)$, while for
i-like cardinals, $\k$ we define
$U(\kappa)$ to be the set of multiples of $\kappa$ in $(\kappa,\ \kappa^+)$.
We define $E$ to be the class of ordinals, $\alpha$, such that letting
$\kappa = card\ \alpha,\ \alpha \in U(\kappa),\ \kappa$ is regular
and either $\kappa$ is inaccessible or ($\kappa$ is s-like
and $\alpha$ is an even multiple of $\k^2$).

For models, $\Cal M,\ Sk_{\Cal M}$
denotes the Skolem hull operator
for $\Cal M$, where the Skolem functions are obtained
in some reasonable fixed fashion.  We often suppress mention of
the membership relation as a relation of a model, but it is usually intended
that it is one.  Thus, $(M,\ A)$
frequently denotes the same model as
$(M,\ \in,\ A)$.

When we have a $\leq^*\text{-increasing}$ sequences of
functions $(\phi_\alpha|\alpha \in X)$, where $X$ is a set of ordinals, and
$\phi$ is a function which is $\leq^*$ one of the $\phi_\alpha$, we let
$scale(\phi)$ denote the least $\alpha \in X$ such that $\phi \leq^*
\phi_\alpha$.
All other notation is introduced as needed (we hope).

%

\pageno=11
\subheading{\S  1.  SINGULAR COMBINATORICS: RESULTS FROM \cite{17}}
\bigskip
\noindent
{\bf (1.1)}\ \ MODEL SEQUENCES AND THE GAMES $G(\th,\ \Cal M,\ p_0)$.
\medskip
Let $\th > \aleph_1$ be regular.  Let $\Cal M = (H_{\nu^+},\ \in,\ \cdots, )$,
where $\nu$ is a singular cardinal, $\nu >> \th$ and $(H_\nu,\ \in)$ models a
sufficiently rich fragment of ZFC.
Let $\si \leq \th$ and let
$(\Cal N_i:i \leq \si)$ be an increasing
continuous elementary tower of elementary substructures of $\Cal M$.
We say that $(\Cal N_i|i \leq \si)$ is
$\bold{(\Cal M,\ \th)\text{{\bf-standard of length }} \si + 1}$
if, letting $N_i \ :=\  |\Cal N_i|$, for all $i \leq \si,
\ card\ N_i = \th,\ \th + 1 \subseteq N_0$, for $i < \si,
\ [N_{i + 1}]^{<\ \th} \subseteq N_{i + 1}$ and, if $i$ is even,
$\Cal N_i \in N_{i + 1}$.

Let $\bold Q$ be a partial ordering (of course, in \S 5, $\bold Q$ will
be $\bold P_\theta$, the upper part of $\bold P$ at $\theta$).  Let
$X$ be dense in $\bold Q$ (in \S 5, below, $X$ will be the dense subclass
of very tidy conditions, see (4.1) and (4.3)).   Let
$\orh {\Cal N} = (\Cal N_i|i \leq \si)$ be $\th\text{-standard}$ with each
$\Cal N_i \prec \Cal M$ (in \S 5, below, $\orh {\Cal N}$ will be
super-coherent (see below)), and let
$q_0 \in Q\ (:=\ |\bold Q|) \cap M\ (:=\ |\Cal M|)$.
The game $G(\th,\ \orh {\Cal N},\ \bold Q,\ X, q_0)$ is defined as follows.\par
{\nar Two players, GOOD and BAD alternate plays.  GOOD plays at
positive even stages (including limit stages); BAD plays at odd stages.
GOOD's moves are conditions, $q_{2i} \in Q \cap M$,
where $0 < i \leq \th$.  For $0 \leq i < \th$, BAD's move at stage $2i + 1$
is a pair,
$(q_{2i + 1},\ \a(i))$, where $q_{2i + 1} \in X,\ q_{2i} \leq
q_{2i + 1},\ \a(i) > sup\ \{\a(j)|j < i\},\
q_{2i},\ q_{2i + 1} \in N_{\a(i)}$.
We require that at all stages $\si \leq \th,\ (q_i:i \leq \si)$ is increasing.
BAD loses if GOOD succeeds in playing $q_\th$.
GOOD loses if at some stage $i \leq \th$, she has no legal move, i.e., there is
no upper bound to the sequence $(q_j:j < i)$.  Of course, this can only occur
if $i$ is a limit ordinal. \bigskip }

We have already hinted at the difficulty
for GOOD at limit stages in the discussion in the
Introduction, preceding the statement of Lemma 5.
\bigskip
\noindent
{\bf (1.2)}\ \  THE SQUARES AND SCALES.
\medskip
\bigskip
\noindent
From \cite{17} (and for singular cardinals of the form
$\aleph_{\tau+\omega}$, using Lemma 3 of the Introduction, above, as
well), we have the following combinatorics for singular cardinals.
\medskip
\noindent
\roster
\item"{(A)}"  A Square on Singular Limits of Limit Cardinals:
\medskip
\noindent
we have a sequence, $(D_\mu: \mu$ is a singular cardinal), where
$D_\mu$ is a club subset of
the singular cardinals below $\mu$ satisfying:
\itemitem {(1)}  $o.t.\ D_\mu < min\ D_\mu$,
\itemitem {(2)}  if $\lambda$ is a limit point of $D_\mu,\
D_\lambda = D_\mu \cap \lambda$.
\itemitem {(3)}  if $\lambda \in D_\mu$ is not a limit point
of $D_\mu$ then $\lambda$ is not a limit of limit cardinals.

\itemitem {(4)}  suppose that $\lambda \in D_{\k_i},\ i = 1, 2$, and let
$j_i$ be such that $\lambda$ is the $j_i^{th}$ member of $D_{\k_i}$.  Then,
$j_1 = j_2$.
\endroster
\medskip
\noindent
If $\tau$ is not a successor ordinal, and
$\k = \aleph_{\tau+\omega}$, conventionally, we let $D_\k\ :=\
\{\aleph_{\tau + n}|n \text{ is odd, } \tau + n \neq 1 \}$.
For such $\k$ we set $\Delta_\k\ :=\ D_\k$.  For $\k$ which are singular
limits of limit cardinals we set $\Delta_\k\ :=\
\bigcup\{\{\lambda\} \cup D_\lambda|\lambda \in D_\k\}$.
\medskip
\roster

\item"{(B)}"  Squares on $(U(\k))^\prime \cap \k^+$,
where $\kappa$ is a singular cardinal
\medskip
for each such $\kappa$, we have a sequence
$(C_\alpha|\alpha \in (U(\k))^\prime \cap \k^+)$
such that each $C_\alpha$ is a club
subset of the set of even multiples of $\k^2$ below $\a$,
of order type less than $\k$, and such that if
$\be \in C_\a$ but is not a limit point of $C_\a$, then $\be$ is not
a limit point of $U(\k)$, and with the usual coherence property:
if $\be$ is a limit point of $C_\alpha,\ C_\beta = C_\alpha \cap \beta$.
\medskip
\item"{(C)}"  Scales on $(\kappa,\ \kappa^+)$,
where $\kappa$ is a singular cardinal):
\medskip
for such $\k$, we have a sequence
$(f^*_\alpha: \alpha \in U(\k))$, where $dom\ f^*_\alpha = D_\k$,
for $\la \in D_\k,\ f^*_\a(\la)$ is an even multiple of $\la^2$
and:
\itemitem {(1)}  if $\kappa < \alpha < \beta,\ \a,\ \be \in U(\k)$
then $f^*_\alpha <^* f^*_\beta$, i.e., for some $\lambda_0 < \kappa$,
whenever $\lambda \in D_\kappa \setminus \la_0,\
f^*_\alpha(\lambda) < f^*_\beta(\lambda)$; further,
if $\alpha \in C_\beta$, then the preceding holds
for {\bf all } $\lambda \in D_\kappa$,
\itemitem {(2)}  whenever $g$ is a function with $dom\ g = D_\k$
and for all $\la \in D_\k,\ g(\la) < \la^+$, for some $\a \in U(\k),\
g <^* f^*_\a$,
\itemitem {(3)}  if $\k$ is a singular limit of limit cardinals,
$\la \in D_\k,\ \a \in U(\k),\ \a^\prime = f^*_\a(\la)$
and $\la^\prime \in D_\k \cap \la$, then $f^*_\a(\la^\prime)
= f^*_{\a^\prime}(\la^\prime)$, and if $\k$ is not a limit of
limit cardinals and $\a,\ \be \in U(\k),\ \la \in D_\k$ and
$f^*_\a(\la) = f^*_\be(\la)$, then $f^*_\a|\la = f^*_\be|\la$,
\itemitem {(4)}  for limit points, $\a$, of $U(\k)$, and
$\la \in D_\k$, $\Phi(\a,\la)\ :=\ \{f^*_\be(\la)|\be \in C_\a\}$
is a final segment of $C_{f^*_\a(\la)}$; further,
on a tail of $D_\k,\ \Phi(\a,\ \la)$ has limit order type.
\endroster
\medskip
Regarding (3), the property given in the second clause follows from the
property given in the first.  Unfortunately, we needed two different
clauses, since we do not have any $f^*_\a$ where
$card\ \a$ is a successor cardinal.
However, the property of the second clause of (3) in fact allows us to
define these according to the following convention.  Once this is done,
in virtue of this definition, we will have the property of the first clause
of (3) even for $\k$ which are not limits of limit cardinals:

{\narrower\smallskip\noindent suppose that
$\la = \aleph_\tau$, where $\tau > 1$ is odd.  Let $\k =
\aleph_{\tau+\omega}$.  Suppose that $\a^\prime = f^*_\a(\la)$ for
some $\a \in U(\k)$.  For $\la^\prime \in D_\k \cap \la$, we
{\it define} $f^*_{\a^\prime}(\la^\prime)$ to be $f^*_\a(\la^\prime)$.
By the second clause of (3), this does not depend on our choice of $\a$.
\smallskip}

Property (4) is the crucial condensation coherence property.  It plays an
important role in the proof, in \cite{17}, of the existence of
super-coherent sequences.
We state this in (1.4), below, as an additional property
of the above combinatorial system, (A) - (C).  Strictly speaking, we never
appeal directly to (4), only to the property of (1.4),
but we do appeal to the following more obvious consequence of (4):

\medskip
\noindent
\ \ \ \ \ \ \ \ \ \ ($4^-$)  on a tail of $D_\k,\
\Phi$ is cofinal in $f^*_\a(\la)$.
\medskip

We close by stating the decodability property of the above system.
As usual, $A$ is as given by Lemma 3 of the Introduction.
\bigskip
\roster
\item"{(D)}"  Decodability of (A) - (C):  for all singular $\k,\ D_\k$
and the systems $(C_\a|\a < \k^+\text{ is a limit point of } U(\k)),\
(f^*_\a|\a \in U(\k))$ are
canonically definable in $L[A \cap \k]$.
\endroster
\medskip
The decodability property is an easy consequence of the fact that
the systems of (B) and (C)
are rather simple modifications of systems which are
canonically constructed in $L$, for singular limits of limit cardinals,
and for $\k$ which are not limits of limit cardinals,
in $L[A \cap \k]$, while the
system of (A) is a simple modification, also given in
\cite{17}, of a constructible system.
\bigskip
\noindent
{\bf (1.3)}\ \ COHERENCE AND SUPERCOHERENCE.
\medskip
Let $\theta > \aleph_1$ be regular. Let
$\nu > cf\ \nu >> \theta$ be such that $(H_\nu,\in) \models$ a sufficiently
rich fragment of ZFC.  Let $\Cal M = (H_{\nu^+},\in,\cdots)$.
Suppose that $\Cal N \prec \Cal M$, where,
letting $N \ :=\  |\Cal N|,
\ card\ N = \th$, and let $\k$ be a cardinal with $\th \leq \k,\ \k \in
N$.  Let $\chi_{\Cal N}(\k) = sup(N \cap \kint)$.

Recall that an Easton set of ordinals is one which is bounded below any
inaccessible cardinal.
For such $\Cal N$ and singular cardinals, $\k$, with
$\theta < \kappa \leq \nu$,
we say that $\kappa$ is $\bold {\Cal N-\text{{\bf{controlled}}}}$ if
there is an Easton set $d$ with $\kappa \in d \in N$.
The Easton sets we have in mind are those consisting of the sets
of cardinals mentioned in some condition in $N$.

We define $p\chi_{\Cal N}$, an analogue of $\chi_{\Cal N}$, defined
on all singular cardinals, $\kappa$, which are $\Cal N-\text{controlled}$.
The definition makes sense for all cardinals $\k \in [\th,\ \nu]$, but we will
only use it for the singulars which are $\Cal N-\text{controlled}$.
If $\kappa \in N$, then of course $\kappa$ is $\Cal N-\text{controlled}$ and
in this case,
$p\chi_{\Cal N}(\kappa)\ :=\ \chi_{\Cal N}(\kappa)$.
Otherwise, $p\chi_{\Cal N}(\kappa)\ :=\ sup\ (\kappa^+ \cap Sk_{\Cal
M}(\{\kappa\} \cup N)).$

The reason that we only consider controlled $\k$ is that one of the results
of \cite{17} gives an alternative characterization of
$p\chi_\N(\k)$ which is central in proving the main result about
the existence of supercoherent sequences (see below).  The alternative
characterization is equivalent only for controlled $\k$.
The restriction to such $\kappa$ is benign, for our
purposes,  since it allows us to handle any cardinal mentioned in any
condition in $\Cal N$.  This is the essential use,
alluded to in the Introduction, just prior to the
\lq\lq DISCUSSION" section, of the fact that the
set of cardinals mentioned in any condition is an Easton set.


Now, let $\modseq$ be $(\Cal M,\ \th)\text{-standard}$ of length $\th + 1$.
For $i \leq \th$, let $\chi_i = \chi_{\Cal N_i},\ p\chi_i =
p\chi_{\Cal N_i}$.  Let $\Cal
N = \Cal N_\th = \bigcup\{\Cal N_i: i < \th\}$, and let $\chi = \chi_\th,
\ p\chi = p\chi_\th$, so
$dom\ \chi = \bigcup\{dom\ \chi_i: i < \th\}$, and for $\k \in dom\ \chi,
\ \chi_\k = sup\ \{\chi_i(\k): \k \in N_i\}$.  Also, for singular
cardinals, $\k \in [\th,\ \nu]$,
which are $\Cal N\text{-controlled},\ p\chi(\k) =
sup\ \{p\chi_i(\k): i < \th\ \&\ \k\text{ is } N_i\text{-controlled}\}$.

Let $\k$ be a singular cardinal, $\k \in dom\ \chi$.
Note that since $cf\ \th = \th > \omega$, there is a club
$D \subseteq \th$ such that for all $i \in D,\ \chi_i(\k) \in C_{\chi(\k)}$.
This motivates the following.
\medskip
\proclaim{Definition}\ \ Let $\M,\ \th$ be as above,
and let $\modseq$ be $(\Cal M,\ \th)\text{-standard}$ of length $\th + 1$.
Let $\Cal N = \Cal N_\th$.
Let $N,\ N_i,\ \chi,\ p\chi,\ \chi_i,\ p\chi_i$ be as above.

Let $\k \geq \th$ be a singular cardinal, $\k \in N$.  $(\N_i: i \leq \th)$ is
$\bold {\M}\text{{\bf{-coherent} at }}\bold \k$
iff for all limit ordinals $\de \leq \th$
with $\k \in N_\de$, there is a club $D \subseteq \de$ such that
for all $i \in D,\ \chi_i(\k) \in C_{\chi_\de(\k)}$.
$(\N_i: i \leq \th)$ is
$\bold {\M}\text{{\bf{-coherent}}}$ if for all
singular cardinals $\k \in N \setminus \th,
\ (\N_i: i \leq \si)$ is $\M\text{-coherent at }\k$.
$(\N_i: i \leq \th)$ is {\bf strongly }
$\bold {\M}\text{{\bf{-coherent}}}$ iff for all
$i < \th$ and all singular cardinals $\k \in N_i,\ \chi_i(\k)
\in C_{\chi(\k)}$.  Finally,
$(\N_i: i \leq \th)$ is {\bf super }
$\bold {\M}\text{{\bf{-coherent}}}$ iff
$(\N_i: i \leq \th)$ is strongly $\M\text{-coherent}$ and for all
limit ordinals, $\si \leq \th$ and all
singular cardinals, $\k$ which are $\Cal N_\si\text{-controlled}$,
for sufficiently
large $i < \si,\ p\chi_i(\k) \in C_{p\chi_\si(\k)}$.
\endproclaim
\medskip
\noindent
{\bf (1.4)}\ \ THE EXISTENCE OF SUPER-COHERENT SEQUENCES.
\medskip
Here is the statement of the main result of
\cite{17} which is the crucial additional property of
the combinatorial system of (1.2).

\proclaim{Lemma}\ \ Let $\theta,\ \nu,\ \Cal M$ be as in (1.3).  Let
$C \subseteq [H_{\nu^+}]^\theta$ be club.  There there is super
$\Cal M\text{-coherent } (\Cal N_i|i \leq \theta)$ with each
$|\Cal N_i| \in C$.
\endproclaim
\bigskip
\noindent
{\bf (1.5)}\ \ AN ADDITIONAL RESULT ABOUT THE SCALES.
\medskip
The following small combinatorial result concerning the scales of
(1.2) is also proved in \cite{17} and will be quite useful in
(4.3), below.

\proclaim{Proposition}  Let $\theta > \aleph_1$ and let $\nu,\ \Cal M$ be as
in (1.4).  Let $d \subseteq [\theta,\nu)$ be an Easton set of cardinals, and
let $\ga$ be a function with domain $d$ such that for all $\k \in d,\
\ga(\k) < \k^+$.  Then, there is a function $\ga^*$ with domain $d$ such that
for all s-like $\k \in d,\ \ga^*(\k) > \ga(\k)$ and such that for all
singular $\k \in d$, letting $\a = \ga^*(\k),\ f^*_\a =^*\ \ga^*|D_\k$.
Further, if $\Cal N \prec \Cal M$ with $(\theta + 1) \cup \{\ga\} \subseteq
|\Cal N|$, then $\ga^* \in |\Cal N|$.
\endproclaim

%

\subheading{\S  2.  PRELIMINARIES ABOUT CONDITIONS}
\bigskip
\noindent
{\bf (2.1)}  CODING AREAS FOR $\eta \in \kint$.
\medskip
We recapitulate, here, some of what we have done in Corollary
4 of the Introduction and (1.2), and
provide some insight into how the coding will work.  First, suppose
that $\eta,\ \k$ fall under one of the following cases.
\medskip
\roster
\item\ \ $\k \geq \aleph_2$ is a successor cardinal, $\eta \in \kint$,
\item\ \ $\k$ is inaccessible, $\eta \in U(\k)$,
\item\ \ $\k$ is a singular cardinal, $\eta \in U(\k)$.
\endroster
Then, by Corollary 4 of the Introduction, for cases (1), (2),
and by (1.2), for case (3),
we have associated to $\eta$ an unbounded subset $b_\eta \subseteq \k$.
In cases (1) and (2), this is the coding area for $\eta$.  In case (3),
it is the main coding area for $\eta$, but we also have
one, and sometimes two auxiliary coding areas for $\eta$ as well,
see below.
\medskip
\noindent
{\bf (2.1.1)} \ \
In case (2), we shall need an additional property of the $b_\eta$.  So, let
$\Cal U\ :=\ \bigcup\{U(\k)|\k \text{ is inaccessible}\}$.  We say that the
system $(b_\eta|\eta \in \Cal U)$ is {\bf tree-like } iff whenever
$\eta_1,\ \eta_2 \in \Cal U$, if $\xi \in b_{\eta_1} \cap b_{\eta_2}$,
then $b_{\eta_1} \cap \xi = b_{\eta_2} \cap \xi$.  In \cite{17} we also prove
the rather simple observation that without loss of generality, we can
assume that $(b_\eta|\eta \in \Cal U)$ is tree-like and has the following
additional property:  $b_\eta = range\ g_\eta$, where $g_\eta$ is a funtion,
$dom\ g_\eta = \{\aleph_\tau|\aleph_tau < card\ \eta \&
\aleph_\tau \text{ is an i-like successor cardinal}
\}$; further, for all $\xi \in b_\eta,\ \xi$ is  a multiple of 4 but
not of 8.

In case (1), if $\k = \mu^+$, it is easy to see that we can,
without loss of generality, assume that the $b_\eta$ have the following
additional properties:  $b_\eta \cap \mu = \emptyset$ and the
members of $b_\eta$ are even ordinals but not multiples of
4; further, if $\mu$ is s-like, then the members of $b_\eta$
are never of the form $\a + 2$, where $\a \in E$.
\medskip
\noindent
{\bf (2.1.2)} \ \
In case (3), (1.2) already gives us $b_\eta$ which have the following
properties.  Once again, the $b_\eta$ are ranges of functions,
$f^*_\eta$, with domain $D_\k$.
Case (3) subdivides according to
whether $\k$ is a limit of singular
cardinals, or of the form $\aleph_{\tau+\omega}$.  In the first
subcase, $D_\k$ is a club subset of singular cardinals below $\k$, whose
order-type is less than its least element.  In the second subcase,
$D_\k$ is the set of $\aleph_\tau$ such that $\tau > 1$ is odd and
such that $\k = \aleph_{\tau+\omega}$.  In both subcases of case (3),
the $f^*_\eta(\la)$ are even multiples of $\la^2$, i.e., they are of
the form $\la^2\iota$ where $\iota > 0$ is even.
\medskip
\noindent
{\bf (2.1.3)} \ \
Recall that a cardinal $\k$ is \bf s-like \rm if $\k$ is singular or
$\k = \aleph_\tau$, where $\tau > 1,\ \tau$ is odd.  If $\k$ is s-like,
$\eta \in U(\k)$,
we set $s_\eta\ :=\
$ the set of multiples of $\k$ in $(\eta,\ \eta + \k^2)$; $s_\eta$ is an
auxiliary coding area for $\eta$ discussed in
in the Introduction, above, and at greater length in (2.4),
(3.4) (E) and (3.5), below.  Finally, if $\k$ is singular and
$\eta \in (U(\k))^\prime$,
we have an additional auxiliary coding area for $\eta$, namely
$C_\eta$, from (1.2) (and so also, implicitly, all of the $s_\a$ for
$\a \in C_\eta$).  This will be used for detecting deactivation of
$\eta$ in a way that is rather important for the limit case of
GOOD's winning strategy in the games of (1.1), and for
determining $\si^{p,\eta}$.  This will also be
discussed more fully in (2.4), below.  It should be noted that unlike the
previous coding areas which are essentially unique to $\eta$, this last is
not, since if $\a \in (C_\eta)^\prime$ then this coding area for
$\a$ is an initial segment of this coding area for $\eta$, and if
$\eta \in (C_\a)^\prime$ then this coding area for $\eta$ is a subset of this
coding area for $\a$.

It is worth recalling
that if $\k$ is a limit of singular cardinals and $\la \in D_\k$ is not a
limit point of $D_\k$, then $\la$ is not a limit of singular cardinals, while
if $\la$ is a limit point of $D_\k$ then $D_\la = D_\k \cap \la$.  It is
also worth recalling (3) of (1.2), and the related convention whereby
we regard $f^*_\nu$ as defined when $\nu = f^*_\a(\la),\
\a \in U((card\ \nu)^{+\omega})$,and $card\ \nu$ is s-like but regular.
\medskip
\noindent
{\bf (2.1.4)} \ \
For regular $\k$ and $\eta$ as in (1) or (2), above, $b_\eta$ will be used
for coding $\eta$ as follows.  If, on a tail of $b_\eta$, we read value 1,
then we will decode value 1 for $\eta$.  Any condition will mention at most
a bounded subset of $b_\eta$, so we will guarantee a tail of 1's on $b_\eta$
by making \lq\lq promises" of the form $(\eta,\ \xi)$, where $\xi < \k$.
Such a pair is a promise to have value 1 at all members of $b_\eta$ above
$\xi$.  By a density argument, (6.1) (2), except for $\k$ is inaccessible
$\eta$ in a bounded subset of $U(\k)$
(the bound is $\be^p(\k)$, see (2.2), (2.3)),
we will have made such a promise {\it whenever} $\eta$ gets value 1.

If $\k$ is a successor cardinal, and we do not have
value 1 on a tail of $b_\eta$, then we will have value 0 on an unbounded
subset of $b_\eta$.
This is also by a density argument, (6.1) (5).
We then decode value 0 for $\eta$.
This will essentially be the procedure when $\k$ is inaccessible,
again, except
for a bounded subset of $U(\k)$.  The situation regarding
the $\eta$ in the bounded set will be discussed more fully in (2.3), and
(3.3), (3.4).
\medskip
\noindent
{\bf (2.1.5)} \ \
For singular cardinals $\k$, and $\eta$ as in case (3), the situation is more
complicated.  Here, any condition which mentions $\eta$ will mention a tail of
$b_\eta$.  In the simplest situation, we will have an $i \in \{0,\ 1\}$
and a tail of $b_\eta$ on which we have value $i$.  It is natural to expect
that when this occurs, we will decode value $i$ for $\eta$.  However, it could
still occur that $\eta$ is deactivated, and that we therefore decode
value !, or sometimes value ? for $\eta$.  We discuss this in (2.4)
and (3.3) - (3.5).  If there is no such tail, then
either we will decode value ? or value ! for $\eta$.  Again this will be
discussed more fully in (2.3), (3.3) - (3.5).
\bigskip
\noindent
{\bf (2.2)}\ \ PROTOCONDITIONS.
\medskip
We define $P(0)$, the class of \lq\lq protoconditions.

\proclaim{Definition}  $p \in P(0)$ iff $p =
(g,\ \be,\ \Xi) = (g^p,\ \be^p,\ \Xi(p))$, and (2.2.1) - (2.2.5), below, hold;
$g$ is the \lq\lq main component", the approximation to the
class function $G$ which are seeking to add (and code down to
a subset of $\aleph_3)$.
\endproclaim
\medskip
\noindent
{\bf (2.2.1)} \ \ There is an Easton set,
$d = d^p$, of cardinals $\geq \aleph_2$, and a function,
$\de = \de^p$ with $dom\ \de = d$, such that for all $\k \in d,\
\k < \de(\k) < \k^+$ and we will have
$g:dom\ g \lra \charset$, with $dom\ g =
\bigcup\{(\k,\ \de(\k))|\k \in d\}$.  $d$ will have the following additional
property:  for singular cardinals $\k \in d$,
there is a tail of $D_\k \subseteq d$.
In addition to the usual characters, 0, 1, we have
the strong deactivator, !, and the weak deactivator, ?, whose roles will
be discussed in (2.3), (2.4), (3.4) and (3.5).
\medskip
\noindent
{\bf (2.2.2)} \ \ $g(\a) \in \{0,\ 1\}$ \bf unless \rm $\k$ is singular
and $\a \in U(\k)$.  However, we have a convention
for systematically abuse of notation for certain $\a$.
\medskip
\noindent
{\bf (2.2.3)} \ \
Recall that $\a \in E$ iff, letting $\k = card\ \a,\
\k > \aleph_1$ is regular, $\a \in U(\k)$ and either $\k$ is inaccessible or
$\k$ is s-like and $\alpha$ is an even multiple of $\k^2$.
In either case, for $i = 0,\ 1$, we take \lq\lq $g(\a) = i$"
as an abbreviation for $(g(\a),\ g(\a + 2)) = (0,\ i)$,
and we take \lq\lq$g(\a) =\ ?$", as an abbreviation for
$(g(\a),\ g(\a + 2)) = (1,\ 0)$.  In the second case {\it only}, we take
\lq\lq $g(\a) =\ !$" as
an abbreviation for $(g(\a),\ g(\a + 2)) = (1,\ 1)$; thus it is our intent
that we can have \lq\lq $g(\a) =\ ?$" for {\it any} $\a \in E$, but
that we have \lq\lq $g(\a) =\ !$" only for those $\a \in E$ whose
cardinalities are s-like.
\medskip
\noindent
{\bf (2.2.4)} \ \
$\be$ is a function with
$dom\ \be = d$, such that for $\k \in d,\ \k < \be(\k) \leq \de(\k)$.
For successor cardinals, $\k \in d,\ \be(\k) = \k + 1$.
The role of $\be(\k)$ for limit cardinals will be made clearer in
(2.3) when we discuss \lq\lq contamination".  For now, we will just say
that $\be(\k)$ is a bound on the contamination in $(\k,\ \k^+)$, not only
in $p$, but in all stronger conditions, $q$.  For cardinals,
$\k \in d$ which are either s-like or inaccessible, we
will have that $\de(\k), \be(\k) \in U(\k)$, and if $\k$ is
inaccessible, we will have that $\be(\k) \geq \k^2$.
\medskip
\noindent
{\bf (2.2.5)} \ \ Finally,
$\Xi$ is the system of \lq\lq promises" which we discussed in (2.1),
above.  $\Xi$ is a set of ordered pairs $(\a,\xi)$ such that
$g(\a) = 1,\ card\ \a$ is regular, $\aleph_2 < \xi < card\ \a$ and
if $card\ \a$ is
inaccessible then ($\a \geq \be(card\ \a)$ and $\a \in U(card\ \a))$.
We shall also require that if $card\ \a
= \la^+$ then $\xi > \la$.
Let $W(p) = dom\ \Xi(p)$.
We let $R(p)\ :=\ \bigcup\{b_\alpha \setminus \xi)|(\alpha,\xi) \in \Xi(p)\}$.
We then require
$g(\zeta) = 1$ for all $\zeta \in R(p)$;
thus $(\a,\ \xi) \in \Xi$ is the \lq\lq promise" to put all 1's
in $b_\a$ from $\xi$ on.
\medskip
\noindent
{\bf (2.2.6)} \ \
In \S 3 we will build to the definition of $P$, by imposing additional
restrictions on the protoconditions.
If $\th > \aleph_2,\ \th$ is regular, then we shall define $P_\th$ in
(4.4).  It is only slightly inaccurate and not at all misleading, at this
point, to say that the main idea is that $d \cap \th = \emptyset$.  The real
point is that $\bold {P_\th}$ is the class of conditions for coding down to a
subset of $\th^+$.

The partial ordering of protoconditions is defined in the most obvious way:
$p \leq q$ iff $g^p \subseteq g^q,\ \be^p \subseteq \be^q$, and $\Xi(p)
\subseteq \Xi(q)$.  This is identical to the definition of the partial
ordering of conditions, in \S 3.
\bigskip
\noindent
{\bf (2.3)}\ \ $X_\gamma$, \lq\lq CONTAMINATION", $\be(\k)$ AND THE
DEACTIVATOR, ?.
\medskip
In (2.1), no coding areas were defined for $\ga$ such that
$\k = card\ \ga$ is a limit cardinal, and $\ga$ is not a multiple of
$\k$.  Strictly speaking, for singular $\k$ and
$\a$ which are multiples of $\k$ but which are
not in $U(\k)$, there was also no coding area defined, but,
except for the multiples of $\k \in [\k,\ \k^2)$, these ordinals are
in $s_\eta$, where $\eta$ is the largest member of $U(\k)$ below $\alpha$.
The multiples of $\k$ in $[\k,\ \k^2)$ are simply ignored.

This is because such ordinals, $\ga$, are not coded directly.
Instead, each such $\ga$ has a set, $X_\ga$, of \lq\lq surrogates", for
coding $\ga$.  Each of the surrogates, $\a$, will have a coding
area $b_\a$ associated with it.
$X_\ga$ will have size $\k^+$, so we have many \lq\lq tries"
at coding $\ga$ correctly.
When $\k$ is singular, this is not entirely
unexpected, since possibly some of the surrogates have been deactivated with
the strong deactivator, !, as in the discussion in the Introduction leading
up to Lemma 5.  Here we discuss the weak deactivator, ?,
and the phenomenon of \lq\lq contamination" which is one of the contexts
in which it arises.  The reasons for calling ? weak and ! strong are discussed
at the beginning of (2.4) and in (2.3.5), below, where we also
discuss another context in which ? arises, for singular
$\k$.  Before doing this, we present the $X_\ga$.
\medskip
\noindent
{\bf (2.3.1)} \ \
For limit cardinals, $\k$, and $\a \in \kint$ which are not multiples of
$\k$, we have sets, $X_\a \in [U(\k)]^{\k^+}$.
If $\k$ is
singular, the $\xi \in X_\a$ are all odd multiples of $\k^2$, i.e., of the form
$\xi = \k^2\iota$, where $\iota$ is odd.  If $\k$ is singular, the system
$(X_\a: \a \in \kint,\ \a \not\equiv 0\ (mod\ \k)) \in L$, while if $\k$ is
inaccessible, each $X_\a = \tilde X_\a \cap \k^+$, where the $\tilde X_\a$ are
classes of ordinals and the relation
$\lq\lq \xi \in \tilde X_\a"$ is canonically $\Sigma_1$ definable over $L$.
When $\k$ is inaccessible, we shall also require the following
property of the $X_\ga$:
\medskip

{\nar(*)\ \ if $\k < \gamma < \zeta< \k^+$ and $\zeta$ is a cardinal in $L$,
then $\zeta = sup(X_\ga \cap \zeta)$.\smallskip}
Finally, for inaccessible $\k$, we take the
$X_\ga$ to partition $U(\k)$, while
if $\k$ is singular, we take the $X_\ga$ to partition the set of odd
multiples of $\k^2$ in $(\k,\ \k^+)$.
\medskip
\noindent
{\bf (2.3.2)} \ \
\lq\lq Contamination" is most easily understood
in the context of inaccessible $\k$.  For such $\k$ and $\a \in U(\k)$, it
can occur that for some inaccessible $\k^\prime > \k$ and
some $\a^\prime \in U(\k^\prime),\
b_\a \cap b_{\a^\prime}$ is unbounded in $\k$.
Further, it could also occur that in some condition $p,\
g^p(\a^\prime) = 1$, and in fact that
$(\a^\prime,\ \xi) \in \Xi(p)$ for some
$\xi < \k$.  This will either prevent us from
having $g^p(\a) = 0$ or from coding this correctly.
When, for other reasons, we are {\it required} to have $g^p(\a) = 0,\
\a$ is said to be \lq\lq contaminated" (by $\a^\prime$) in $p$.
We cannot prevent such contamination, but we will
define conditions in such a way (see (3.2) (A),
below) that \par
{\nar ($\ast$):  fewer than $\k$ many
$\a \in \kint$ are contaminated. \medskip}
\noindent
Typically, contamination occurs here because $\k$ was added to
$d^p$ after the promise $(\a^\prime,\ \xi)$ had
already been made.
\medskip
\noindent
{\bf (2.3.3)} \ \
When $\k$ is singular, we shall also have the phenomenon of
contaminated ordinals.  It may
occur, for singular $\k$, and conditions, $p$, with $\k \in d^p$, that
for some $\xi \in U(\k) \cap \de^p(\k)$, one of the
following holds:
\medskip
\roster
\item"{(1)}"  There are $x_1 \neq x_2$ and $Y_1,\ Y_2$, cofinal subsets of
$b_\xi$, such that for $i = 1,\ 2$ and $\z \in Y_i,\ g^p(\z) = x_i$,
\medskip
\item"{(2)}"  $i \in \{0,\ 1\}$, and for other reasons we are required
to have $g^p(\xi) = i$,
but on a cofinal set of $\z \in b_\xi,\ g^p(\xi) = 1 - i$  (see (3.3) below,
the definition of \lq\lq forced to be $i$", for
what these \lq\lq other reasons" are).
\endroster
\medskip
\noindent
If $g^q(\xi) =\ !$, then $\xi$ will be deactivated anyway.  If, however, this
fails, then $\xi$ is contaminated in $q$.
Once again, contamination will occur only for a bounded set of $\xi$, though
here this is a simple observation which does not require a special property of
the conditions, as in the inaccessible case.
Here again, typically, contamination arises due to the fact that
an unbounded set of information below $\k$ was part of a condition before
$\k$ was mentioned.

In both the inaccessible and the singular case,
$\be^p(\k)$ is the sup of
the contaminated ordinals $\xi \in \kint$ (see (3.4) (C)).
Once $\be^p(\k)$ has been specified in a condition, no
further contamination is allowed in any stronger condition, $q$, since
$\be^q(\k) = \be^p(\k)$.  In both the inaccessible and the
singular case, we require, (3.4) (B), that if $\a$ is contaminated
then $g^p(\a) =\ ?$.
\medskip
\noindent
{\bf (2.3.4)} \ \
We can now specify how the $\a \in X_\ga$ are used to code $\ga$.
If $\ga \in dom\ g^p$, then we will have $g^p(\ga) \in \{0,\ 1\}$.
In the inaccessible case, we will have that if $\a \in X_\ga \setminus
\be^p(\k)$, then $g^p(\a) = g^p(\ga)$ (see (3.4) (A), and one clause
of the definition of \lq\lq forced to be $i$").  In the singular case,
things are somewhat more complicated, since even for $\a \in X_\ga
\setminus \be^p(\k)$, we can have $g^p(\a) \in \{?,\ !\}$.  However,
as part of the definition of condition ((3.4) (A), again), we will have that
for such $\a,\ g^p(\a) \bold \neq 1 - g^p(\ga)$.  We will show, by a density
argument, (6.1) (4), that when $\Cal I$ is generic and
$G$ is the union of the $g^p$ for $p \in \Cal I$, there will be a cofinal
set of $\a \in X_\ga$ such that $G(\a) = G(\ga)$.  Thus, in decoding,
there is a common definition:  decode for $\ga$ the unique value
$i \in \{0,\ 1\}$ such that we have value $i$ on a cofinal subset
of $X_\ga$.
\medskip
\noindent
{\bf (2.3.5)} \ \
To conclude, we should mention the other way the weak deactivator,
?, can occur.  For singular $\k$, in addition to  occurring at contaminated
$\a$, it can occur at other $\a \in U(\k)$, but only if on a tail of
$b_\a$ the value ? occurs.  The reason for this has exactly to do with the
density argument we just mentioned.  As will become clearer
in (2.4), the strong deactivator at $\a$ can contribute to deactivating
larger ordinals.  This is not the case for the weak deactivator, ?.
Thus, the weak deactivator, ?, can play the role of \lq\lq safe, neutral
filler", and does not present the \lq\lq potential danger" of forcing us to
deactivate ordinals we want to preserve
as \lq\lq active", to get value in $\{0,\ 1\}$, such as the
cofinally many $\a \in X_\ga$ we need for the preceding.
\bigskip
\noindent
{\bf (2.4)}\ \ THE STRONG DEACTIVATOR, !, AND THE GENERIC SCALES.
\medskip
Suppose that $\k$ is singular and $\a \in U(\k)$.
We have already mentioned most of the elements of this discussion:
\medskip
\roster
\item\ \ when $g^p(\a) =\ !$, not only $\a$, but also
all members of $U(\k)$ in the open interval
$(\a,\ h^p(\a))$ are deactivated (if $g^p(\a) =\ !$ and $h^p(\a) = \a,\ \a$
itself is still deactivated);
this is one of the senses in which ! is the strong deactivator.
\item\ \ $h^p(\a) = scale(\si^{p,\a})$.
\item\ \ if $g^p(\a) =\ !$ this can contribute to making $g^p(\nu) =\ !$, for
certain larger $\nu \in U(\k)$; this is the other sense in which ! is the
strong deactivator.
\item\ \ when $\a$ is a limit point of $U(\k),\ C_\a$ is an additional,
auxiliary coding area used for detecting strong deactivation.
\item\ \ detecting strong deactivation and
lengths of deactivated intervals (by (2),
above, this amounts to the same thing as decoding the $\si^{p,\a}$) is one of
our major preoccupations.
\endroster
\medskip
\noindent
{\bf (2.4.1)} \ \
We now put these elements together and lay the groundwork for (3.5), omitting,
for now, some of the finer points related to certain $\nu \geq \de^p(\k)$ for
which $\si^{p,\nu}$ will nevertheless be defined.  We should say, at the
outset, that $\si^{p,\a}$ will be defined whether or not we end up
having $g^p(\a) =\ !$, but that this is just for convenience, since the only
case in which it has any significance is when this occurs; when $g^p(\a)
\neq\ !$, we ignore $\si^{p,\a}$ and take $h^p(\a)$ to be $\a$.
As we have already mentioned, we are grateful to the
referee for emphasizing the point of view that
the $\si^{p,\a}$ are really potential members of generic scales
which we are forcing as we do the coding.
In almost all cases,
we will have $\si^{p,\a} \geq^* f^*_\a$; the exceptions are
discussed in (3.5), (3.6).
\medskip
\noindent
{\bf (2.4.2)} \ \
We will have two other functions,
$\upsilon^{p,\a}$ and $\pi^{p,\a}$ and that,
in most cases, for $\la \in D_\k \cap d^p$, we take
$\si^{p,\a}(\la)\ :=\ max(\upsilon^{p,\a}(\la),\ \pi^{p,\a}(\la))$.  Looking at
$\upsilon^{p,\a}$ amounts to considering what happens \lq\lq from below", on
$b_\a$.  Looking at $\pi^{p,\a}$ amounts to considering what happens
\lq\lq to the left", on $C_\a$.  These are two of the ways in which $\a$
could be strongly
deactivated, and are two of the places we have to look to detect
strong deactivation.
\medskip
\noindent
{\bf (2.4.3)} \ \
Before developing this, however, there is a third way
in which $\a$ can be strongly
deactivated, and we deal with this first, since it is
simplest, and directly related to (1), above.  $\a$ is
$p\text{-interval-strongly-deactivated}$
if it is in
a deactivated interval, $(\nu,\ h^p(\nu))$, for some $\be^p(\k) \leq \nu <
\a$.  When
this occurs, we take $\nu$ least possible and set $\si^{p,\a}\ :=\
\si^{p,\nu}$, without considering the $\upsilon^{p,\a},
\ \pi^{p,\a}$.  Thus, when
$\a$ is $p\text{-interval-stronly-deactivated}$,
\lq\lq only this counts", even if it turns out
that it is also deactivated in one of the two other ways we now discuss.
In terms of Lemma 4.3, this corresponds to the $\a$ between $\de$ and
the $t^p_2(\k)$ of (4.2), and, roughly speaking, to limit stages of GOOD's
winning strategy.
\medskip
\noindent
{\bf (2.4.4)} \ \
$\a$ is $p\text{-strongly-deactivated}$ on $b_\a$ (\lq\lq from below")
iff on a tail of $\xi \in b_\a,\ g^p(\xi) =\ !$.
Typically, this occurs when $\a$ {\it didn't have to be deactivated}, when
we are strongly
deactivating $\a$ intentionally, to be sure that we are
able to strongly
deactivate other, larger, $\a$ which will be more problematical,
see the discusssion of $(\ast)$ in the Introduction.  This corresponds to
some of the work in (4.3) (beyond the $t^p_2(\k)$, of (4.3)) and all of the
work of (6.1), and, roughly speaking, to successor stages in GOOD's winning
strategy.
\medskip
\noindent
{\bf (2.4.5)} \ \
$\a$ is $p\text{-strongly-deactivated}$ on
$C_\a$ iff it is a limit point of $U(\k)$ and
on a cofinal subset of $\nu \in C_\a,\ g^p(\nu) =\ !$.  Typically, this
occurs in situations where we {\it really needed to deactivate} $\a$
and we are happy to find that we prepared for this by
strongly deactivating enough
members of $C_\a$.  This corresponds to the portion of the work
in (4.3) dealing with $\de^p(\k)$ and to the situation of $\a =
\de(\k)$ in (4.5), and roughly speaking, to limit stages in
GOOD's winning strategy.
\medskip
\noindent
{\bf (2.4.6)} \ \
It remains only to give the main idea of the
definitions of the $\upsilon^{p,\a}$ and the $\pi^{p,\a}$
(there are some fine points which can be deferred until
the official definition in (3.5)).  The main idea for the
$\upsilon^{p,\a}(\la)$ is that this should be $h^p(f^*_\a(\la))$.
The fine points arise when $f^*_\a(\la) \geq \de^p(\la)$.
The main idea for the $\pi^{p,\a}(\la)$ is that this should be
$sup\ \{\si^{p,\nu}(\la)|\nu \in C_\a\}$.  The fine points arise because
we want this sup to be $\geq f^*_\a(\la)$, but $\leq \de^p(\la)$.
\bigskip
\noindent
{\bf (2.5)}\ \ OVERVIEW OF THE DECODING PROCEDURE.
\medskip
Let $\h = \h_A$ be the (class) characteristic function of $A$.  Our forcing
will produce a generic class function $G$ with domain $\subseteq OR$
and range $\subseteq \charset$.
\noindent
We will code $A$ into $G$ on odd ordinals, i.e., we
shall have that for non-successor ordinals $\de$ and $n < \omega,\  G(\de + 2n
+ 1) = \h(\de + n)$.\par
Of course, we want to recover $G$ from $G|\aleph_3$ by \underbar {decoding}.
This is done by recursion on cardinals, $\k$.
The basic recursion step is to go from $G|\k$ to $G|(\k,\ \k^+)$,
when $\k \in CARD$.  This will involve a nested recursion across $(\k,\ \k^+)$.
The procedure for obtaining $G|(\k,\ \k^+)$ from $G|\k$ will be
uniform within each of the following classes of cardinals:  inaccessibles,
singulars, and successors.
Thus, at limit cardinals, $\m$, we can piece together $G|\m$ from the $G|\k,\
\k < \m$, and continue.  The recursion step for successor
cardinals is provided by (2.1.4).  As noted there,
for inaccessibles, this
also {\it essentially} gives the way we obtain $G_0$,
which we now discuss.\par

For limit cardinals, $\k$, it will simplify matters if, in decoding
$G|(\k,\ \k^+)$, we
have available not only $G|\k$, but also an auxiliary function,
$G_0$, which represents the first stage in defining $G|(\k,\ \k^+)$.
The role of
$G_0$ can best be understood by discussing the broad outline of how we finally
obtain $G|(\k,\ \k^+)$.
For inaccessible cardinals, this is a\lq\lq two-pass" process.
For singular cardinals, it is a \lq\lq three-pass" process.

The first pass involves decoding the information provided by
$G|\k$ on $b_\eta$ without regard to
the analogous information for the $\nu \in (\k,\ \eta)$.
$G_0$ represents the outcome of this \lq\lq first
pass".  For inaccessibles, even this first pass involves a recursion,
since we have to decode the $b_\eta$ as we go.  For singulars, however,
there is no recursion involved in the first pass, but there definitely is a
recursion involved in the second pass for singulars, where we deal with
the strong deactivator, !, and the generic scales.  The second pass
for inaccessibles and the third pass for singulars are analogous, in that this
is where we deal with contamination, and define $G$ on
the non-multiples of $\k \in \kint$.

%

\subheading{\S 3.  THE CODING CONDITIONS:  DEFINITIONS}
\bigskip

We build to the definition of the class of coding
conditions, ${\bold P}$, in (3.5) - (3.6).
In our original treatment we had stronger properties,
which appear below as (4.1) (A) and (B+), in place of (3.5) (C) and (D).
The latter are technical weakenings of the properties of (4.1),
which are designed to allow us to prove, in (4.3), that the very tidy
conditions, those with the properties of (4.1) are dense.
\bigskip
\noindent
{\bf (3.1)}
\medskip
We recall some terminology and conventions from the Introduction
and (2.2).
A cardinal $\k$ is {\bf s-like} if it is singular
or of the form $\aleph_\tau$, with $\tau > 1$ and odd.  Next,
let $\kappa$ be a regular uncountable
cardinal and let $\alpha \in \kint$.  Recall that $\alpha \in E$
if $\a \in U(\k)$ and either $\kappa$ is inaccessible or
$\kappa$ is s-like.  Formally, for conditions $p$
and $\a \in E$ we shall have $g^p(\a) \in \{0,1\}$, but
recall the convention from (2.2.3) involving the use of $\a + 2$ as an
\lq\lq extra bit" for $\a \in E$.
Naturally, we have taken care not to
assign any other \lq\lq coding duties" to the $\a + 2$ where $\a \in E$.
\bigskip
\noindent
\proclaim{(3.2)\ \ Definition}
\medskip
Suppose $p = (g,\be ,\Xi) = (g^p,\be^p,\Xi(p)) \in P(0)$, where $P(0)$ is
as in (2.2).

$p \in {\tilde P}$ iff the following
properties, (A) and (B)  are satisfied.\endproclaim
\medskip
\roster
\item"{(A)}"  For all regular $\kappa^\prime$,
there are fewer than $\k^\prime$ many $\a$ such that for some
$\xi  < \k,\ (\a,\ \xi) \in \Xi(p)$ (note:
$\kappa^\prime$
need not be a member of $d$),
\medskip
\item"{(B)}" If $\a \in dom\ g^p,\ \k = card\ \a$ is singular and $\a$
is a multiple of $\k^2$:
\medskip
\ \ \ \ (1)  $f^*_\a <^* \de^p|D_\k$,

\ \ \ \ (2)  if $g^p(\a) \in \{0,1\}$, then on a tail of $b_\a,\ g^p(\zeta)
= g^p(\a)$,

\ \ \ \ (3)  if $\be^p(\k) \leq \a$ and $g^p(\a)\ =\ ?$, then on a tail of
$b_\a,\ g^p(\zeta)\ =\ ?$.
\endroster
\bigskip
\noindent
{\bf (3.3)}
\medskip
Suppose $p \in {\tilde P}$.  First, consider $\alpha$ such that
$\kappa = card\ \alpha$ is a limit cardinal, and
suppose that $\alpha \in X_\gamma$.
We say that $g^p(\alpha)$ {\bf is forced to be 0 (resp. 1)} if
$g^p(\gamma) = 0$ (resp. 1).
We also say that
$g^p(\alpha)$ {\bf is forced to be 1} if for some $\gamma \in
R(p),\ \alpha \in X_\gamma$.
Finally, drop the restriction on $\kappa$.
If $\alpha = 2\alpha^\prime + 1$,
then we say that
$g^p(\alpha)$ {\bf is forced to be 0 (resp. 1)} if $\alpha^\prime
\not\in A\text{ (resp. } \alpha^\prime \in A)$.

If $\k = card\ \a$ is inaccessible and $\a$
is a multiple of $\k$, then $\a$ is {\bf contaminated by} $\bold\tau$ if
$\tau \in W(p),\ (\tau,\ \xi) \in \Xi(p)$ for
some $\xi < \k,\ b_\tau \cap b_\a$ is cofinal in $\k$
and $g^p(\a)$ is forced to be 0; $\a$ is {\bf contaminated } iff for
some $\tau$ it is contaminated by $\tau$.  Because the system of $b_\a$
is tree-like for $\a \in \Cal U$ (see (2.1.1)), it is easy to see that
any $\tau \in W(p)$ contaminates at most one $\a \in
(\k,\k^+)$.  Therefore, (3.2) (A) gives that there are fewer than
$\k$ many $\a \in (\k,\k^+)$ which are contaminated.

If $\k = card\ \a$ is singular, $d^p \cap D_\k$ is cofinal in $\k$, then
$\a$ is {\bf contaminated } iff $\a$ is a multiple of
$\k^2,\ g^p(\a) \neq\ !$ and one of the
following holds:
\medskip

\ \ \ \ (1)  there are $x_1 \neq x_2$ and cofinal subsets $Y_1,\ Y_2
\subseteq b_\a \cap dom\ g^p$ such that for $\xi \in Y_i,\ g^p(\xi) = x_i$,
\medskip
\ \ \ \ (2)  $g^p(\a)$ is forced to be 0 (resp. 1) but on a cofinal
subset of $b_\a \cap dom\ g^p,\ g^p(\zeta) = 1$ (resp. 0).

Here, it is easy to see that at most $\k$ many $\a \in \kint$
are contaminated, since if $\a > scale(\de^p|D_\k)$ then $\a$
cannot be contaminated.  Also, note that $\a$ which are contaminated because
of (2) are {\it odd} multiples of $\k^2$ since they are members of some
$X_\ga$.
\bigskip
\noindent
\proclaim{(3.4) Definition}
\medskip
If $p
\in {\tilde P}$, then $p \in P^*$ iff the following properties (A) - (E)
hold.\endproclaim
\medskip
\roster
\item"{(A)}"  If $g^p(\a)$ is forced to be $i$ and $g^p(\a) \in \{0,1\}$
then $g^p(\a) = i$,
\medskip
\item"{(B)}"  If $\a \in dom\ g^p$ and $\a$ is contaminated, then
$g^p(\a) =\ ?$,
\medskip
\item"{(C)}"  For limit $\k \in d,\ \be^p(\k) =
sup\ \{\a \in (\k,\k^+)|\a \text{ is contaminated}\}$,
\medskip
\item"{(D)}"  If $\k \in d^p,\ \k$ is inaccessible, we also define $\be^p_1(\k)
\ :=\ sup\ \{\a \in (\k,\k^2)|\a \text{ is contaminated}\}$ and
we require:

\ \ \ \ (1)  $g^p(\beta^p_1 + \k\om) = 0,\ g^p(\beta^p_1 + \k\si) = 1$, for
all limit ordinals $\si < \k$,

\ \ \ \ (2)  $g^p|(\beta^p_1 + \kappa\omega,\ \kappa^2)$
codes a well-ordering of $\kappa$ in type $\beta^p(\kappa)$
on odd successor multiples of $\k$ in
$(\be^p_1(\k) + \k\omega,\ \k^2)$, and codes
\newline\ \ \ \ \ \ \ \ $A \cap \be^p(\k)$ on
even successor
multiples of $\k$ in $(\be^p_1(\k) + \k\omega,\ \k^2)$.
\medskip

\item"{(E)}"  Suppose $\kappa \in d$ is s-like.
Suppose that $\be^p(\k) \leq \a < \de^p(\k)$, with
$\a$ a multiple of $\k^2$.
Let $\Ga(\ ,\ )$ denote the
G\"odel pairing function.  Let $H^p(\a)\ :=\
\{(\xi,\ \zeta) \in \kappa\times\kappa|g^p(\a + \k (1+\Ga (\xi,\
\zeta))) = 1\}$.  We then require that $H^p(\a)$ is a well-ordering of a
subset of $\k$ which lies in $L[A \cap \k]$, and,
if $\k$ is singular, we also require that it is the $<_{L[A \cap
\k]}\text{-least}$ well-ordering of a subset of $\k$ in its order type.

We let $h^p(\a)\ :=\ $ the least multiple of $\k^2 \geq$ the order type of
$H^p(\a)$.
We further require:
\medskip
\ \ \ \ (1)  $h^p(\a) \geq \a$,

\ \ \ \ (2)  if $\k$ is singular, $\be^p(\k) \leq \be < \a$,
and $\be$ is a multiple of $\k^2$ then
$h^p(\a) \geq h^p(\be)$, and if $h^p(\be) \geq \a$, then $h^p(\a) = h^p(\be)$,

\ \ \ \ (3)  If $g^p(\a)\ \neq\ !$, then we require that $h^p(\a)$ has
the smallest value consistent with (1) and (2) (which, it will be easy to
see, from (3.5) (A), will be $\a$).
\endroster
\bigskip

\noindent
{\bf (3.5)}
\medskip
\proclaim{(3.5.1) Definition}
Fix $p \in P^*$ and singular $\kappa \in d^p$.
Suppose that $\a \in U(\k) \setminus \be^p(\k)$.  We are mainly interested
in the case where $\a \leq \de^p(\k)$, but it will be useful to have
the definition in the more general context.  This results in somewhat more
complicated definitions; we will also give the simpler definitions that
result when we restrict to $\a \leq \de^p(\k)$.

Let $g = g^p,\ h = h^p,\
f = f^*_\a,\ \beta = \beta^p(\k)$.  We first define some additional functions,
$\upsilon^{p,\a},\ \pi^{p,\a},\ \sigma^{p,\a}$, with domain
$D_\kappa \cap d^p$.

First, for
$\lambda \in dom\ \upsilon^{p,\a}$, if $f(\lambda) \in dom\ g$, we set
$\upsilon^{p,\a}(\lambda) = h^p(f(\lambda))$; otherwise,
$\upsilon^{p,\a}(\lambda)
= \de^p(\lambda)$.  Note that if $\a \in dom\ g$, then on a tail
of $\la \in D_\k \cap d,\ \upsilon^{p,\a}(\la) = h^p(f(\la)) \geq f(\la)$.
Thus, if it is {\bf not} the case that $f \leq^* \upsilon^{p,\a}$ then
$\a \geq \de^p(\k)$ and there is no tail of $b_\a \subseteq dom\ g$.

We now define $\pi^{p,\a},\ \sigma^{p,\a}$ by simultaneous recursion
on $\a$; at the same time we define three properties,
$Pr^p_1,\ Pr^p_2,\ Pr^p_3$, by defining, by recursion on $\a$
when $Pr^p_i(\a)$ holds.  We use $Pr^p(\a)$ as an abbreviation for
$Pr^p_1(\a)$ or $Pr^p_2(\a)$ or $Pr^p_3(\a)$.  We say that
$\a$ is $\bold{p\text{{\bf -interval-strongly-deactivated\ }}}$ iff
$\a \leq \de^p(\k)$ and $Pr^p_1(\a)$ holds.
We say that $\a$ is $\bold{p\text{{\bf -strongly-deactivated}}}$
{\bf on} $\bold{b_\a}$
iff $\a \leq \de^p(\k)$ and $Pr^p_2(\a)$ holds.  Finally,
we say that $\a$ is
$\bold{p\text{{\bf -strongly-deactivated}}}$ {\bf on } $\bold{C_\a}$
iff $\a \leq \de^p(\k)$ and $Pr^p_3(\a)$ holds.
We say that $\a$ is
$\bold p\text{{\bf -strongly-deactivated}}$ iff it is
$p\text{-strongly-deactivated}$
on $b_\a$ or it is $p\text{-interval-strongly-deactivated }$
or it is $p\text{-strongly-deactivated on }C_\a$.  Thus, $\a$ is
$p\text{-strongly deactivated iff } \a \leq \de^p(\k)$ and $Pr^p(\a)$ holds.

We turn, now, to the recursive
definition of the two above-mentioned functions, and the three properties.
$Pr^p_1(\a)$ holds just in case
there is $\nu \in U(\kappa) \cap [\beta,\ \a)$, such
that $Pr^p(\nu)$ holds and
$scale(\sigma^{p,\nu}) > \a$.
If $Pr^p_1(\a)$ holds, let $\nu$ be the least
witness to this.  In this case, we set $\pi^{p,\a} = \pi^{p,\nu}$, and
$\sigma^{p,\a}\ :=\ \sigma^{p,\nu}$.

Thus, for the definition of the two functions,
we can assume $Pr^p_1(\a)$ fails.  In this case,
for $\lambda \in dom\ \pi^{p,\a},\ \pi^{p,\a}_1(\lambda)\ :=\
min(\de^p(\la),f(\lambda)),\ \pi^{p,\a}_2(\la)\ :=\
sup\ \{\sigma^{p,\nu}(\lambda)|\nu \in C_\a\})$
and $\pi^{p,\a}(\la)\ :=\ max(\pi^{p,\a}_1(\la),\pi^{p,\a}_2(\la))$.
Finally, for $\lambda \in dom\ \sigma^{p,\a},\
\sigma^{p,\a}(\lambda) = max(\upsilon^{p,\a}(\lambda),\pi^{p,\a}(\lambda))$.

We conclude by defining when the other two properties,
$Pr^p_2(\a),\ Pr^p_3(\a)$ hold.  $Pr^p_2(\a)$ holds
iff on a tail of $\lambda \in D_\kappa,\ Pr^p(f(\lambda))$ holds.
$Pr^p_3(\a)$ holds iff $\a$ is
a limit of multiples of $\kappa^2$, and $Z_\a$ is cofinal in $\a$, where
$Z_\a = \{\nu \in C_\a|Pr^p(\nu)\text{ holds}\}$.

Of course, more than one of these may be true for $\a$.  However, if
$Pr^p_1(\a)$, \lq\lq only this counts", in
terms of how $\sigma^{p,\a}$ is defined.
It is also possible that, letting $\de = \de^p(\k),\
\de)$ is $p\text{-deactivated}$.
This is clear in the case of interval deactivation
and deactivation on $C_\de$.  Deactivation on
$b_\de$ is only possible if a tail of $b_\a \subseteq dom\ g$.

The simplifications which arise when we restrict to $\a \leq \de^p(\k)$
are mainly that we can remove the definitions of $Pr^p_2(\k)$ and
$Pr^p_3(\a)$ from the recursion which gives us the definitions of
$\pi^{p,\a}$ and $\si^{p,\a}$, by changing the definition of
$Pr^p_2(\a)$ to be:  \lq\lq on a tail of $\lambda \in D_\kappa,\
f(\lambda) \in dom\ g\ \&\ g(f(\lambda)) =\ !$", and for $Pr^p_3(\a)$,
by changing the definition of $Z_\a$ to:  \lq\lq $\{\nu \in C_\a|g(\nu) =\
!\}$."  The reasons will be clear from (A), below.  We can also drop from the
definition of $Pr^p_1(\a)$ the requirement that $Pr^p(\nu)$ holds.
\endproclaim

A disquieting possibility is that $Pr^p(\a)$ holds for {\bf all } $\a \in
U(\k) \setminus \be^p(\k)$.
In Remark 2 of (4.3) we shall show that this cannot occur.
We are now ready for the definition of $P$.

\proclaim{(3.5.2) Definition}  $p \in P$ iff $p \in P^*$,
property (D), below, holds, and whenever $\kappa,\ \a$,
etc., are as above, and $\de = \de^p(\k)$,
the following properties (A) - (C) hold:
\endproclaim
\medskip
\roster
\item"{(A)}"  If $\a < \de$, then $g(\a)\ =\ !$ iff $\a$ is
$p\text{-strongly-deactivated}$,
\item"{(B)}"  If $\a < \de$,
and $g(\a) =\ !$, then $h^p(\a) = scale(\sigma^{p,\a})$,
\item"{(C)}"  If $\lnot(\de^p|D_\k \leq^* f^*_\de)$, then $\de$ is
$p\text{-strongly-deactivated},\ \de^p|D_\k \leq^* \sigma^{p,\de}$;
further, letting $\ga = scale(\sigma^{p,\de})$, whenever
$\eta \in U(\k)$ with
$\de < \eta < \ga$, if $b_\eta \cap dom\ g^p$ is cofinal
in $\k$, then on a tail of $\xi \in b_\eta \cap dom\ g^p,\ h^p(\xi) \leq
f^*_\ga(card\ \xi)$,
\item"{(D)}"  If $\la \in d^p$ is s-like, the following
set has power $\leq \la$:

$$\{\si^{p,\a}(\la)|\si^{p,\a}(\la) \text{ is defined}\}.$$

\endroster
\bigskip
\proclaim{Remark}
The substantive part of (D) concerns those
$\si^{p,\a}(\la)$ which are $> \de^p(\la)$.
As indicated at the beginning of this section, our original definition
of $P$ required that all the $\si^{p,\a}(\la) \leq \de^p(\la)$ and
that, with the notation of (C), above, $\de^p|D_\k =^* f^*_\de$.
Instead, we have opted to relax this requirement and show,
in \S 4, that these properties hold on a dense set.  With this in mind,
(D) is clearly a necessary condition to be able to extend $p$ to a
condition with these properties.  (C) is a technical property, formulated
with the same aim.
\endproclaim
\medskip
\proclaim{(3.6)  Definition}
\medskip
If $p,\ q \in P$, we set $p \leq q$ iff
$g^p \subseteq g^q,\ \be^p
\subseteq \be^q,\
\Xi(p) \subseteq \Xi(q)$.
\endproclaim
\proclaim{Remark 1}  Note that if $p \leq q$ then $h^p \subseteq h^q$.
Note, also, that if $\a \geq \be^p(\k),\
\a \in dom\ g^p \cap U(\k)$ then $\upsilon^{p,\a} =^*
\upsilon^{q,\a},\ \pi^{p,\a}_i =^* \pi^{q,\a}_i,\ i = 1,2$,
and therefore $\pi^{p,\a} =^* \pi^{q,\a}$ and $\si^{p,\a} =^*
\si^{q,\a}$.  It is also easy to see that if $\de = \de^p(\k)$,
then $\pi^{p,\de}_2 =^* \pi^{q,\de}_2$.  It is possible that
$\upsilon^{q,\de}(\la) > \upsilon^{p,\de}(\la)$; this will occur exactly when
$\de^p(\la) < f^*_\de(\la) < \de^q(\la)$.  Similarly,
$\pi^{p,\de}_1(\la) < \pi^{q,\de}_1(\la)$
just in case $\de^p(\la) < f^*_\de(\la) < \de^q(\la)$.
Thus, we could
have $\upsilon^{q,\de}(\la) >
\upsilon^{p,\de}(\la)$ on a tail of $\la$.
It is also clear
that $\de$ is $p\text{-interval-strongly-deactivated }$
just in case it is $q\text{-interval-strongly-deactivated}$, and similarly for
deactivation on $C_\de$.  However, it is possible that $\de$ is
$q\text{-strongly-deactivated on } b_\de$
without being $p\text{-strongly-deactivated on } b_\de$.

The situation is similar for $\a \in U(\k) \setminus \de + 1$.  It is easy to
see that if $Pr^p_i(\a)$ holds then $Pr^q_i(\a)$ holds, and that
$\upsilon^{p,\a} \leq^* \upsilon^{q,\a},\ \pi^{p,\a}_i \leq^* \pi^{q,\a}_i$
and therefore that $\pi^{p,\a} \leq^* \pi^{q,\a},\ \si^{p,\a} \leq^*
\si^{q,\a}$.
\endproclaim
\proclaim{Remark 2}
In virtue of (3.5)(B), above, letting $\de = \de^p(\k)$,
we define $h^p(\de)$ by $h^p(\de)\ :=\ scale(\si^{p,\de})$.
\endproclaim
%
\proclaim{Remark 3}  Let $\de = \de^p(\k),\ \ga = scale(\de^p|D_\k)$, and
suppose that $\a \in U(\k) \cap \ga$. Note that this occurs exactly when
$\a \in U(\k),\ \de < \a$ and
$b_\a \cap dom\ g^p$ is cofinal in $\k$.  Thus,
for such $\a$ there is already an unbounded set of information
imposed by $p$ on $b_\a$, which might require us to deactivate
$\a$, and the question arises of how far this deactivation should go.
However, if this occurs, then we have the hypotheses of (3.5)(C), above,
and so $\de$ is $p\text{-strongly-deactivated}$.  Further, since
(3.5)(C) gives us that $\de^p|D_\k \leq^* \sigma^{p,\de}$, and
$\a < \ga,\ \a < scale(\sigma^{p,\de})$.  Thus, $Pr^p_1(\a)$ holds.
Suppose, now that $q \geq p$ and $\de < \de^q(\k)$.  By Remark 1, above,
$\si^{p,\de} \leq^* \si^{q,\de}$.  Thus, in such
$q \geq p,\ \a$ will already be
$q\text{-interval-strongly-deactivated}$ by $\de$.  Finally, since
$\a < scale(\de^p|D_k)$ and, by hypothesis, $b_\a \cap dom\ g^p$ is cofinal in
$\k$, on a tail of $\xi \in b_\a \cap dom\ g^p, h^p(\xi) \leq
f^*_\ga(card\ \xi)$.  Thus, as far as such $\a$ are concerned, $\de$ already
provides the essentials of the deactivation information.
\endproclaim

%

\subheading{\S 4.  THE CODING CONDITIONS:  BASIC LEMMAS}
\bigskip
\noindent
\proclaim{(4.1)  Definition}
\medskip
If $p \in P,\ p\ \text{is {\bf tidy} }$
iff for all s-like $\k \in d^p$, (A), below, holds and for all
singular $\k \in d^p$, (B), below holds.
\medskip
\roster
\item"{(A)}"  if $\a \in U(\k)$ with $\be^p(\k) \leq \a
< \de^p(\k)$, then $h^p(\a) \leq \de^p(\k)$,
\item"{(B)}"  $\de^p|D_\k \leq^* f^*_\de$, where $\de = \de^p(\k)$.
\endroster
\medskip
If $p \in P$, then $p$ is {\bf very tidy } iff
for all s-like $\k \in d^p$, (A), above, holds
and for all singular
$\k \in d^p$, (B+), below, holds.
\medskip
\roster
\item"{(B+)}" $\de^p|D_\k =^* f^*_\de$, where $\de = \de^p(\k)$.
\endroster
\endproclaim
\medskip
\proclaim{Remark 1}  Suppose that $p$ is tidy.  We argue that for
all singular $\k \in d^p$, all $\a \in U(\k) \cap [\be^p(\k),\ \de^p(k)$
and all $\la \in D_\k \cap d^p,
\ \si^{p,\a}(\la) \leq \de^p(\la)$.  It suffices,
of course, to prove this for the $h^{p,\a}$ and the $\pi^{p,\a}$.  For
the $h^{p,\a}$,
if $f^*_\a(\la) \notin dom\ g$, then $h^{p,\a}(\la) = \de^p(\la)$,
so suppose that $f^*_\a(\la) \in dom\ g$.  Then, $h^{p,\a}(\la) =
h^p(f^*_\la(\a))$, and by (A) above (with $\la$ in place of $\k$ and
$f^*_\a(\la)$ in place of $\a$), the latter is $\leq \de^p(\la)$, as
required.  For the $\pi^{p,\a}$, we work by induction on $\a$, with the
induction hypothesis being the statement of the remark, i.e., the statement
for the $\si^{p,\nu}$, with $\nu < \a$.  But then, the conclusion is
immediate by the definition of $\pi^{p,\a}$:  clearly,
$\pi^{p,\a}_1(\la) \leq \de^p(\la)$,
and $\pi^{p,\a}_2(\la)$
is the $sup$ of things all $\leq \de^p(\la)$ and so the conclusion is
clear.
\endproclaim
\medskip
\proclaim{Remark 2}  If $p$ is very tidy then for singular$\k \in d^p$,
it is easy to see that, with the convention of Remark 2 of (3.6),
$h^p(\de^p(\k)) \leq \de^p(\k)$.  This is clear from Remark 1
and the fact (which is just a restatement
of (B+)) that $\de^p(\k) = scale(\de^p|D_\k)$.
\endproclaim
\medskip
\proclaim{Remark 3}  If $p$ is very tidy, $p \leq q$ and for all s-like,
regular $\la \in d^p,\ h^q(\de^p(\la)) = \de^p(\la)$ then for all
s-like $\k \in d^p, h^q(\de^p(\k)) = \de^p(\k)$.  This is easily argued by
induction on the rank of $\k$ in the well-founded relation
\lq\lq $\la \in D_\k$" .  The basis is the hypothesis.
Let $\eta\ :=\ \de^p(\k)$.  By the induction hypothesis, we have that
$h^{q,\eta} =^* \de^p|D_\k$.  Clearly $\pi^{q,\eta}_2 =^* \pi^{p,\eta}_2$,
and by Remark 2, $\pi^{p,\eta}_2 =^* \de^p|D_\k$.  Clearly,
$\pi^{q,\eta}_1 =^* \de^p|D_\k$, and the conclusion is then immediate.
\endproclaim

\noindent
{\bf (4.2)}
\medskip
The following material will be  helpful in
both (4.3) and \S 5.  If $p \in P$, and $t$ is a function with
$d^p \subseteq dom\ t$, we say that $t {\bold\ covers\ } p$
iff whenever $\k \in d^p$ is singular,
$\a \in U(\k) \cap [\be^p(\k),\ \de^p(\k)]$
and $\la \in D(\k) \cap d^p,\ \si^{p,\a}(\la)
< t(\la)$.  If $q \in P$, we we say that $q {\bold\ covers\ } p$
iff $p \leq q$ and $\de^q$ covers $p$.  If $q \in P$ and $t$ is a
function with $dom\ t = d^q$, we say that $q {\bold\ dominates\ }
t$ iff for all s-like $\la \in d^q$
we have $t(\la) < \de^q(\la)$.

Next, we define still more functions associated with a $p \in P$.
For singular $\k \in d^p$, and $\la \in D_\k \cap d^p$,
we let $t^p_{\k,1}(\la)\ :=\
sup\ \{\si^{p,\a}(\la)|\a \in U(\k) \cap [\be^p(\k),\ \de^p(\k)\}$.
We note that by (3.5)(D), $t^p_{\k,1}(\la) < \la^+$.
For regular, s-like $\k \in d^p$, we let $t^p_1(\k)\ :=\
sup\ \{h^p(\a)|\a \in U(\k)\ \&\ \a < \de^p(\k)\}$.
Again, by (3.5) (D), for regular, s-like $\k \in d^p,\ t^p_1(\k) < \k^+$.
For singular $\k \in d^p$, we let $t^p_1(\k)\ :=\
scale(t^p_{\k,1})$.  Clearly, for singular $\k \in d^p,\ t^p(\k) < \k^+$.
We also define the $t^p_{\k,2}$ and the $t^p_2$ analogously, but based on the
function $\de^p$; thus, $t^p_{\k,2}\ := \de^p|D_\k \cap d^p$,
and $t^p_2(\k)\ :=\ scale(t^p_{\k,2})$, for singular $\k \in d^p$,
while for s-like regular $\k \in d^p,\ t^p_2(\k)\ :=\ \de^p(\k)$.
Note that whenever these functions are defined, we have
$t^p_{\k,2}(\la) \leq t^p_{\k,1}(\la)$ and $t^p_2(\k) \leq t^p_1(\k)$.

Now, let $\theta,\ \Cal M,\ \Cal N$, etc., be as in (1.2), (1.3),
and suppose that $p \in N$.
Then it is obvious that the $t^p_i \in Sk_{\Cal M}(N)$
and that the $t^p_{\k,i} \in Sk_{\Cal M}(N \cup\{\k\})$, and therefore, that

$$(\ast)\text{  for all s-like } \k \in d^p,\ t^p_1(\k) < p\chi_\N (\k),$$

\noindent
since $(\ast)$ holds for {\it any} $t \in Sk_{\Cal M}(N \cup \{\k\})$
in place of $t^p_1$.

\proclaim{(4.3)\ \ Lemma}
\medskip
If $p \in P$, $t$ is a function with
$dom\ t = d^p$ and for all $\k \in d^p,\ t(\k) < \k^+$, then
there is very tidy $q \in P$ with $p \leq q$ and such that for all
s-like $\k \in d^p,\ t(\k) < \de^q(\k)$.
\endproclaim
\demo{Proof}  We shall prove this in the way that will be most useful for
(5.1).  Choose regular $\theta \geq \aleph_2$,
and let $\Cal M,\ \Cal N$ be as above for this $\theta$.  We have
just observed that since $p \in N,\ p\chi_{\Cal N}$
is everywhere $\geq t^p_1$; similarly, since
$t \in N$, if we construct very tidy $q \geq p$ such that

$$(\ast)\text{  for all s-like\ } \k \in d^p, \de^q(\k)  >
p\chi_{\Cal N}(\k),$$

\noindent
then $q$ will be as required.  This is the approach we shall take; we shall
choose such an $\Cal N$, and construct $q$ satisfying $(\ast)$
and such that whenever $\k \in d^p$ is s-like,
$g^q(\de^p(\k)) =\ !$.  Our approach
to this will be to take $\ga = p\chi_\N$, and to let $\ga^*$ be as given by
(1.5) for this $\ga$, and to take $d^q = d^p,\ \Xi(q) = \Xi(p),\ \be^q = \be^p,
\ \de^q(\la) = \ga^*(\la)$, for all s-like $\la \in d^p$ and
$\de^q(\la) = \de^p(\la)$ for all other $\la \in d^p$.  Recall that,
for all singular $\k \in d^p$, letting $\nu = \ga^*(\k),\
\ga^*|D_\k =^* f^*_\nu$.  This makes it clear that we will have
(B+) of (4.1) and that $q$ will satisy $(\ast)$.
Thus, in order to complete the proof, it will suffice
to verify that (A) of (4.1) holds and that $q \in P$,
since it will then be clear that $q \geq p$.

Before going further, it will be useful to exploit (4.1) (B+) further.
Suppose that $u$ is a function with domain $ =^* D_\k$.
We make the following observations:
\medskip

\roster
\item"{(A)}"  Suppose that for all
$\la \in dom\ u,\ u(\la) \leq \ga^*(\la)$.  Then, clearly,
$scale(u) \leq \nu$.

\item"{(B)}"  Suppose further $b$ is a cofinal subset of $D_\k$ and
that for $\la \in b,\ u(\la) = \ga^*(\la)$.  Then, again clearly,
$scale(u) = \nu$.
\endroster
\medskip
An important property of the way we will
define $q$ is:
\medskip
\roster
\item"{(C)}"  $g^q(\eta) =\ !$, for all $\de^p(\k)
\leq \eta < \de^q(\k)$, whether or not $\k$ is singular.
\endroster
\medskip
By (C) and the definition of $t^p_2$
if $\k$ is singular, then,
\medskip

\roster
\item"{(D)}"  for $t^p_2(\k) \leq \eta < \de^q(\k)$,
with $\eta \in U(\k),\ \eta$ will be
$q\text{-strongly-deactivated on } b_\eta$, if for no other reason.
\endroster
\medskip

It remains to see that for singular $\k$ and $\de^p(\k) \leq
\eta < t^p_2(\k)$, with $\eta \in U(\k)$, we will still have
that $\eta$ is $q\text{-strongly-deactivated}$.  It is here that we will appeal
to (3.5) (C).  Let $\de\ :=\ \de^p(\k)$.  If $\de^p|D_\k
\leq^* f^*_\de$, then $t^p_2(\k) = \de$, so there is nothing to verify
in this case.  If, on the other hand, the above fails, then,
by (3.5) (C), $\de$ is $p\text{-strongly-deactivated }$ and $\de^p|D_\k \leq^*
\sigma^{p,\de}$.  Since it is clear that $\sigma^{p,\de} \leq^*
\sigma^{q,\de}$, this means that we will have that $\de$ is
$q\text{-strongly-deactivated}$ and that we will have $h^q(\de) \geq
scale(\de^p|D_\k)$, and $scale(\de^p|D_\k)$ is just $t^p_2(\k)$.
This in turn means that if $\de < \eta < t^p_2(\k)$,
with $\eta \in U(\k)$,
then we will have that $\eta$ is $q\text{-interval-strongly-deactivated, }$
as required.

The preceding guarantees that we can carry out our plan of making (C) hold,
while respecting (3.5).  It remains to complete the definition of
the $g^q|[\de^p(\k),\de^q(\k))$ and
to define the $h^q(\eta)$, for $\de^p(\k) \leq
\eta \leq \de^q(\k)$, with $\eta \in U(\k)$.
This will be done by recursion on the rank of $\k$ in the well-
founded relation:\ \ \lq\lq $\la \in D_\k$", so the basis is when $\k$
is regular.  In view of (C), we must carry out the
following:
\medskip
\roster
\item"{(1)}"  define $h^q(\eta)$ and $g^q|s_\eta$ for $\eta$ as above,
\item"{(2)}"  define $g^q(\xi)$ for $\de^p(\k) < \xi < \de^q(\k)$ which
are not covered by (C) nor by (1).
\endroster
\medskip

As far as (2) is concerned, in all cases, we shall have $g^q(\k)\
:=\ 1$ unless it is forced to be 0.  As far as (1) is concerned,
once we have computed $h^q(\eta)$, (3.4) (E) tells us how to define
$g^q|s_\eta$.  Thus, it remains to compute the $h^q(\eta)$ and verify that
the computed value is consistent with (3.4), (3.5) and (4.1) (A).
It will then be clear that $q$ is a very tidy condition, and of course,
that $p \leq q$, completing the proof of the Lemma.  Of course, (3.5) (B)
tells us that for singular $\k$, in order
to compute the $h^p(\eta)$, it suffices to compute the
$\sigma^{p,\eta}$.  This will be done by recursion on $\eta$, within the
recursion on $\k$.

We now appeal to (A) and (B).  Our induction hypotheses are
\medskip
\roster
\item"{(E)}"  for all relevant $\la < \k$ and
$\de^p(\la) \leq \eta^\prime < \de^q(\la),\ h^q(\eta^\prime)
\leq \ga^*(\la)$,
\item"{(F)}"  for all relevant $\de^p(\k) \leq \eta^\prime < \eta$ and
{\bf all, not just on a tail of, } $\la \in D_\k \cap d^p,\
\sigma^{q,\eta^\prime}(\la) \leq \ga^*(\la)$.
\endroster
\medskip

Now, (E) guarantees that for all $\la \in D_\k \cap d^p,\ h^{q,\eta}(\la)
\leq \ga^*(\la)$.  Similarly, (F) guarantees that for all $\la \in
D_\k \cap d^p,\ \pi^{q,\eta}_2(\la) \leq \ga^*(\la)$, and therefore,
for all such $\la,\ \sigma^{q,\eta}(\la) \leq \ga^*(\la)$.  This preserves
the induction hypothesis, (F), and then, by (A), $scale(\sigma^{q,\eta})
\leq \ga^*(\k)$, which preserves the induction hypothesis, (E), at least
as far as $\eta$.  As indicated above, this completes the proof that
$q \in P$ and $p \leq q$.
\enddemo

\proclaim{Remark 1}
The $q$ we obtain
depends, obviously, on $p$ and on the $\Cal N$ we choose, so, we naturally
denote it as $q(p,\ \Cal N)$.  A very plausible choice for $\Cal N$ is to take
it to be some sort of Skolem hull in $\Cal M$ of $p,\ t$ and possibly some
additional elements, e.g., all the members of $\theta$, for some cardinal
$\theta$.  We return to this in (5.1), (5.2), below.  In this connection, it
is also worth pointing out that if $p \in |\Cal N|,\
\Cal N^* \prec \Cal M$, and, letting $N^*\ :=\ |\Cal N^*|$, if
$[N^*]^{<\theta} \cup N \cup \{\Cal N\} \subseteq |N^*|,\ card\ N^* =
\theta$, then clearly $q \in N^*$.
\endproclaim

\proclaim{Remark 2}  We are now in a position to
show, as promised at the end of (3.5.1), that in no condition
$q \in P$ do we have that for a tail of $\a \in U(\k) \setminus \de^q(\k) + 1,\
Pr^q(\a)$ holds, where $\k \in d^q$ is singular.
What we show, in fact, is that if $p \in P$ is very tidy,
$\k \in d^p$ is singular, $\a \in U(\k) \setminus \de^p(\k) + 1$ then
$Pr^p(\a)$ fails.  This suffices, since as noted in the last sentence of
Remark 1 of (3.6), if $q \leq p$, then for $\a \in U(\k)$, if $Pr^q_i(\a)$
holds then $Pr^p_i(\a)$ holds.

This will also complements Remark 3 of
(3.6).  since if $p \in P$ is very tidy, and $\k \in d^p$ be singular, then
letting $\de = \de^p(\k)$, we have that $\de = scale(\de^p|D_\k)$, so there
are no ordinals $\a$ of the sort dealt with in Remark 3 of (3.6).  We would
like to know that no others share the property pointed out there, of being
$q\text{-strongly-deactivated}$ in {\it any } $q \geq p$ with $\de^q(\k)
\geq \a$.  Actually, this will follow from the extendability properties
developed in (6.1), but showing that if $Pr^p(\a)$ holds then $\a \leq
\de^p(\k)$ will in fact be quite useful for (6.1).

If $\a \in U(\k) \setminus \de + 1$, since $p$ is very tidy,
there is a tail of $b_\a$ disjoint from $dom\ g^p$, so $Pr^p_2(\a)$ fails.
Again, since $p$ is tidy, we cannot have
$h^p(\eta) > \a$, for any $\eta < \de^p(\k)$,
and by Remark 3 of (4.1), we cannot have $h^p(\de) > \a$.
We conclude by induction on $\a$, so suppose that for all $\ga \in (\de,\ \a)
\cap U(\k),\ Pr^p(\ga)$ fails.  Clearly, then $Pr^p_3(\a)$ cannot hold, and
$Pr^p_1(\a)$ cannot be witnessed by any $\ga \in (\de,\ \a)$.  But we have
already argued that $Pr^p_1(\a)$ cannot be witnessed by any $\eta \leq \de$, so
the proof is complete.
\endproclaim
\noi
{\bf (4.4)}\ \ FACTORING.
\medskip
Let $\theta$ be a regular cardinal, $\theta \geq \aleph_2$,
and let $p \in P$.  We set $W_\theta(p)\ :=\ W(p) \setminus \theta^+$,
and we set $\Xi_\theta(p)\ :=\
\{(\a,\ max\ (\xi,\ \theta))|(\a,\ \xi) \in \Xi(p) \& \alpha \in
W_\theta(p)\}$.
We let $W^\theta(p)\ :=\ W(p) \cap \theta^+,\
\Xi^\theta(p)\ :=\ \Xi(p) \cap W^\theta(p) \times \theta$ and we let
$R^\theta(p)\ :=\ \bigcup \{b_\a \cap [\xi,\theta)|(\a,\xi) \in \Xi(p)
\&\ \xi^p(\a) < \theta,\ \theta^+ \leq \a \}$.
We are now ready to define the upper and lower
parts of $\bold P$, relative to $\theta$, which give the Factoring Property
of $\bold P$.
\medskip
\proclaim{Definition}\ \ For $p$ and $\theta$ as above, we set
$(p)_\theta = (g^p|d^p \setminus \theta,\
\be^p|d^p \setminus \theta,\ \Xi_\theta(p))$.
Note that $(p)_\theta \in P$.  We let $\bold {P_\theta} = \{(p)_\theta|p \in
P\}$, with the restriction of $\leq$.  Thus, $\bold {P_\theta}$ is the class
of conditions for coding down to a subset of $\theta^+$.  Note that $\bold P =
\bold {P_{\aleph_2}}$, since for $p \in P,\ (p)_{\aleph_2} = p$.

We also define $(p)^\theta$, for $p \in P$:\ \
$(p)^\theta = (g^p|\theta,\ \be^p|\theta,\ \Xi^\theta(p),\ R^\theta(p))$,
and we let $\dot {P}^\theta\ :=\
\{((p)^\theta,\ (p)_\theta)|p \in P\}$.  Thus, $\dot {P}^\theta$
is a (proper class) $\bold P_\theta$~-~name for a subset of
$\{(p)^\theta|p \in P \}$.  We have guaranteed that the latter is a set
by replacing
$\{(\a,\ \xi)|(\a,\ \xi) \in \Xi(p)\ \&\ \xi < \theta,\ \theta^+ < \a \}$ by
$R^p(\theta)$.  Of course, our intention is to have
$\dot P^\theta$ be the name of the underlying set of a partial subordering
of $(\{(p)^\theta|p \in P\}, S)$, where $(p)^\theta\ S\ (q)^\theta$ iff
$g^p|\theta \subseteq g^q,\ \be^p|\theta \subseteq \be^q,\
R^\theta(p) \subseteq R^\theta(q),\
\Xi^\theta(p) \subseteq \Xi^\theta(q)$.  We let
$\dot {\bold P}^\theta$ be the name for this subordering.
\endproclaim

In fact we can cut $\dot{\bold P}^\theta$ down to a set name, as follows.
Let $n < \omega$
be such that all relevant notions about $\bold P$ are $\Sigma_n$.  Let $\chi
\ >>\ \theta$ be such that $(H_\chi,\ \in)$ reflects all $\Sigma_n$ formulas.
The set name is then simply \newline
$\{((p)^\theta,\ (p)_\theta)|p \in P \cap
H_\chi\}$.  What makes $\dot P$ a name is the linkage between the \lq\lq top"
and the \lq\lq bottom", which is what
guarantees that $\dot {\bold P}^\theta$ will code down to a subset of
$\aleph_3$ the subset $\dot B$ of $\theta^+$ added by $\bold P_\theta$.
The $R^\theta(p)$ is one feature of this linkage.  The following is then
clear:
\proclaim{Lemma}\ \ (FACTORING) $\bold P \cong
\bold P_\theta\ast\dot{\bold P}^\theta$.
\endproclaim
\noindent
\proclaim{(4.5)\ \ Lemma}
\medskip
Let $\theta$ be as in (4.4), and suppose
that $\si \leq \theta$ is a limit ordinal and
$(p_i|i < \sigma)$ is an increasing sequence from $P_\theta$.
Let $d\ :=\ \bigcup\{d^{p_i}|i < \si \},\
g\ :=\ \bigcup\{g^{p_i}|i < \si \},\ \be\ :=\
\bigcup\{\be^{p_i}|i < \si \},\ \Xi\ :=\ \bigcup\{\Xi(p_i)|i < \si \}$.
For $\k \in d$, let $\de(\k)\ :=\
\bigcup\{\de^{p_i}(\k)|i < \si\ \&\ \k \in d^{p_i}\}$.  If $\k
 \in d$ is singular, set $\eta \in Z(\k)$ iff $\eta \in C_{\de(\k)}$
and $g(\eta) =\ !$.  Let $i \in I(\k)$ iff for some $i \leq j < \si$ and some
$\eta \in Z(\k),\ \de^{p_i}(\k) \leq \eta$ and for all
$\la \in D_\k \cap d^{p_i},\
\de^{p_i}(\la) \leq \si^{p_j,\eta}(\la)$.

Suppose, further, that $(p_i|i < \si )$ has
the following properties.
\medskip
\roster
\item  $\{i < \si|p_i\text{ is tidy}\}$ is cofinal.
\item  For all singular $\k \in d,\ I(\k)$
is cofinal in $\si$ (this, of course, implies that $Z(\k)$
is cofinal in $\de(\k)$).
%
\endroster
\medskip

Then, $p\ :=\ (g,\be ,\Xi) \in P_\theta$ and is the least upper bound
for $(p_i|i < \si )$.
\endproclaim
\demo{Proof}  We will concentrate on showing that (3.5) (C) holds.
This is the heart of the matter for verifying that $p \in P_\theta$, as
verifying that the other clauses hold is totally routine, and once
we know that $p \in P_\theta$ it is clear that it is the least upper bound.

So, let $\k$ be as in the statement of the Lemma, and adopt the other
notation there.  Also, let $\a\ :=\ \de(\k)$, and note that $Z(\k)$
is just the $Z_\a$ of (3.5).  As observed in the parenthetical
remark to hypothesis (2) of the statement of the Lemma,
$Z_\a$ is therefore cofinal in $\a$, which means that $\a$ is
$p\text{-strongly-deactivated on } C_\a$.
So, it remains to verify the last clause of (3.5) (C).

For this, we first note that
for all $\la \in D_\k \cap d$, we have

$$ (\ast)\
\ \pi^{p,\a}_2(\la) \geq \de(\la).$$

This is because, since $I(\k)$ is
cofinal in $\si$, if $\la \in D_\k \cap d,\
\de(\la) = sup\ \Delta$, where
$\Delta\ :=\ \{\de(i)(\la)|i \in I(\k)\ \&\ \la \in d^{p_i}\}$.
Now, let $i \in I(\k)$ with $\la \in d^{p_i}$.  Let
$i \leq j < \si$ and $\eta \in Z(\k)
\setminus \de^{p_i}(\k)$ be as guaranteed by the fact that
$i \in I(\k)$.  Then, $\pi^{p,\a}(\la) \geq \si^{p_j,\eta}(\la) \geq
\de^{p_i}(\la)$.

We now complete the proof by verifying the last clause of
(3.5) (C).  So, let $\ga\ :=\ scale(\si^{p,\a})$.  Note that in virtue
of $(\ast)$, and the fact that for all relevant $\la,\ \si^{p,\a}(\la)
\geq \pi^{p,\a}_2(\la)$, we have that $f^*_\ga \geq^* \si^{p,\a} \geq^*
\pi^{p,\a}_2 \geq^* \de|D_\k$.  Now, if $\eta \in U(\k)$
with $\a < \eta < \ga$, if $b_\eta \cap dom\ g^p$ is cofinal in
$\k$ and $\xi \in b_\eta \cap dom\ g^p$, then $\xi < \de(card\ \xi)$.  But
then there is $i < \si$ with $\xi < \de^{p_i}(card\ \xi)$, and in virtue of
(1), we can assume that $p_i$ is tidy.  But then $h^{p_i}(\xi)
< \de^{p_i}(card\ \xi)$, by (4.1) (A), and since $h^{p_i}(\xi)
= h^p(\xi)$, the conclusion is clear.
\enddemo
\proclaim{Remark}  In addition to the hypotheses of the Lemma, suppose that
$\si < \th$, that $\Cal M$ is as in (1.2), (1.3), that
$\Cal N \prec \Cal N^* \prec \Cal M$ are such that, letting $N\ :=\ |\Cal N|,\
N^*\ :=\ |\Cal N^*|$, we have that for all $i < \si,\ p_i \in N$
and that $[N^*]^{<\theta} \cup N \cup \{\Cal N\} \subseteq N^*$.  Then $p \in
N^*$.
\endproclaim
\bigskip
\noindent
{\bf (4.6)}\ \ DECODING.
\medskip
We now complete the sketch of the decoding procedure
given in (2.5)
by supplying the details of the decoding for limit cardinals.
For singulars, in addition to the case of decoding the generic,
we also treat the case of decoding $g|[\k,\ \de^*)$ from a $g$
defined on a large enough domain below $\k$.  We will give the
details of the situation when we encounter it.  This case is needed
in (6.1) (7) (b), below.
We treat the inaccessible case first.
\bigskip
\noindent
{\bf (4.6.1)}
\medskip
Assume that $\k$ is inaccessible, and that we are given
$G$, a function with $range\ G \subseteq \{0,1,?,!\}$ and
$dom\ G = \bigcup\{(\la,\la^+)|\aleph_2 \leq \la < \k,\ \la\text{ a
cardinal}\}$ such that for some generic ideal $\Cal I$ in $\bold P,\
G = \bigcup\{g^p|\k\ |p \in \Cal I\}$.

We will first obtain an ordinal $\be$ which will be the
common value of the $\be^p(\k)$ for $\k \in d^p,\ p \in \Cal I$.
Recall that, by Lemma 3 of the Introduction
and the discussion preceding it, for all $\a \in (\k,\ \k^+)$, if
$\k \leq \nu \leq \a$ and in $L[A \cap \nu],\ card\ \a = \k$, then
we obtain $(b_\eta|\k < \eta \leq \a)$ canonically from $\a$ in
$L[A \cap \nu]$.  In particular, we have, in $L,\ (b_\a:\a < \k^2)$.
So, we first decode $G$ on $U(\k) \cap \k^2$.
For such $\a$, we let $G_0(\a) = 1$ iff on a tail of $\xi \in b_\a,\
G(\xi) = 1$.  Otherwise, we set $G_0(\a) = 0$.  Now, by (3.4) (D) (1),
there is a largest $\nu \in (U(\k))^\prime \cap (\k,\ \k^2)$ such
that $G_0(\nu) = 0$ and, further that this $\nu$ is of the form
$\eta + \k\om$, where $\eta \in \{\k\}\cup U(\k)$.
Also, by (3.4)(D)(2), we have that $G_0|(\eta,\ \k^2)$ codes a well-ordering of
$\k$ on odd successor multiples of $\k$ in this interval.  We take
$\be\ :=$ the order-type of this well-ordering.
By (3.4)(D)(2), again, we have that $\be$ is the common value of
$\be^p(\k)$ for $\k \in d^p,\ p \in \Cal I$ and that $A \cap \be$ is coded
by $G_0$ on even successor multiples of $\k$ in this interval.
\bigskip
\noindent
{\bf (4.6.2)}
\medskip
We can now define $G(\a)$ by recursion on $\a$ for
$\a \in U(\k) \setminus \be$.  We first define $\nu(\a)$ and obtain
$A \cap \nu(\a)$.  We shall have that in $L[A \cap \nu(\a)],\ card\ \a =
\k$, so that we have $b_\a$ available.  If $\a$ is not a cardinal in
$L$, we let $\nu(\a)\ :=\ card^L\ \a$.  Otherwise, we let
$\nu(\a) = \a$.  Recall from (2.3)
that if $\k < \zeta< \k^+$ is a cardinal in $L$,
then for all non-multiples of $\k, \ga$
with $\k < \ga < \zeta,\ X_\ga \cap \zeta$ is cofinal in $\zeta$.
This allows us to define $A \cap \nu(\a)$ as follows.
If $\nu(\a) \leq \be$, then we already have $A \cap \be$.  Otherwise, if
$\k < \xi < \nu(a)$, let $\ga\ :=\ 2\xi + 1$.
Then, $\ga < \nu(\a)$, so $X_\ga \cap (\be,\ \nu(\a)) \neq \emptyset$, and
we have $\xi \in A$ iff for some (all)
$\eta \in X_\ga \cap (\be,\ \nu(\a)),\ G(\eta) = 1$.
Thus, we have $A \cap \nu(\a)$, and therefore, in
$L[A \cap \nu(\a)]$, we have $b_\a$.  We set $G(\a)\ :=\
1$ iff on a tail of $\xi \in b_\a,\ G(\xi) = 1$; otherwise,
$G(\a) = 0$.  This completes the recursion.
\bigskip
\noindent
{\bf (4.6.3)}
\medskip
We can then define $G$ on
$(\k,\ \k^+) \setminus U(\k)$ by $G(\ga)\ :=\ i$ iff
for some (all) $\a \in X_\ga \setminus \be,\ G(\a) = i$.  Finally, we can go
back and define $G$ on $U(\k) \cap \be$ as follows:
let $\ga \not\in U(\k)$ be such that $\a \in X_\ga$.  If $G(\ga) = 1$ and
on a tail of $\xi \in b_\a,\ G(\xi) = 1,\ G(\a)\ :=\ 1$.  If $G(\ga) = 0$
and on a cofinal subset of $\xi \in b_\a,\ G(\xi) = 0$, then
$G(\a)\ :=\ 0$.  Otherwise, we set $G(\a)\ :=\ ?$.
This completes the decoding procedure in the inaccessible case.
\bigskip
\noindent
{\bf (4.6.4)}
\medskip
So, assume next that $\k$ is singular.  We first treat the generic case.
Assume that $G, \Cal I$ are as above.  This time, there are
\lq\lq three passes"
in the definition of $G$.  The first is relatively straightforward:
we ignore deactivated intervals, deactivation on $C_\a$, we ignore
the $\pi^{p,\a}$ and the $\si^{p,\a}$, and just organize
the information \lq\lq from below"
provided by $G$.  This is done simultaneously, for all members of $U(\k)$
at once.  No recursion is involved.
Recall here that we have $A \cap \k$ at our disposal, and that
all of the coding apparatus is present in $L[A \cap \k]$.
The second pass proceeds by recursion on $\a \in U(\k)$,
and takes into account all of the above; we also define
$G|s_\a$.  At the end of the second pass, we will have defined a function
$G_1$ on all multiples of $\k$ in $(\k,\ \k^+)$.
In the third pass, we then go back,
define $G$ on the non-multiples of $\k$, detect
contamination, define $\be$ (which, once again, will be the common value of
the $\be^p(\k)$ for $p \in \Cal I$ with $\k \in d^p$) and revise the definition
of $G_1$ below $\be$.

For the first pass, if $\a \in U(\k)$ and
there is no tail of $b_\a$ on which $G$ is constant, set $G_0(\a)\ :=\
?$; otherwise, if $x \in \{0,\ 1,\ ?,\ !\}$ and $G$ has constant value $x$ on
a tail of $b_\a$, set $G_0(\a)\ :=\ x$.  We also define $\Upsilon^\a$
at this time,
by letting $\Upsilon^\a(\la)\ :=\ $ the order-type of the well-ordering
coded by $G|s_{f^*_\a(\la)}$, for $\la \in D_\k$.
\bigskip
\noindent
{\bf (4.6.5)}
\medskip
For the second pass, we assume that we have defined $G_1$ on all multiples of
$\k$ below $\a$ in such a way that for
$\k < \nu < \a$, with $\nu \in U(\k),\
G_1|s_\nu$ codes a well-ordering of $\k$ in order type $\geq \nu$.
We let $H(\nu)\ :=$ the order type of this well-ordering.
We also assume that we have defined functions $\pi^\nu,\ \si^\nu$ with domain
$D_\k$, for such $\nu$, \lq\lq correctly" (i.e. according to (3.5), and what
now follows), so far. If there is such a
$\nu < \a$ with $H(\nu) > \a$ we
take the least such $\nu$ and
set $G_1(\a)\ :=\ !,\ \si^\a(\la) =
\si^\nu(\la)$, for all $\la \in D_\k$, and we define $G_1|s_\a$ to code the
$<_{L[A \cap \k]}\text{-least}$ well-ordering of $\k$ in type $H(\nu)$.

So, assume there is no such $\nu$.
Let $f\ :=\ f^*_\a$ and define $\pi^\a_1\ :=\ f$ and for $\la \in D_\k$,
define $\pi^\a_2(\la)\ :=\ sup\ \{\si^\nu(\la)|\nu \in C_\a\}$ and
$\pi^\a(\la)\ :=\ max(\pi^\a_1(\la),\ \pi^\a_2(\la)),\ \si^\a(\la)
= max(H^\a(\la),\ \pi^\a(\la))$.
If, on a tail of $\xi \in b_\a,\
G(\xi) =\ !$, then we already had $G_0(\a) =\ !$ and we maintain
$G_1(\a)\ :=\ !$.  However we also set $G_1(\a)\ :=\ !$, if
$\a$ is a limit of multiples of $\k^2$ and $Z_\a$ is cofinal in
$\a$ where $Z_\a\ :=\ \{\nu \in C_\a|G_1(\nu) =\ !\}$.
In all other cases, we maintain $G_1(\a)\ :=\ G_0(\a)$.
If $G_1(\a) =\ !$, we define $G_1|s_\a$ to code the
$<_{L[A \cap \k]}\text{-least}$ well-ordering
of $\k$ in type $scale(\si^\a)$.
In all other cases we define $G_1|s_\a$ to code the
$<_{LA \cap \k]}\text{-least}$ well-ordering of $\k$ in type $\a$.
This completes the
recursive definition of $G_1$ on the multiples of $\k$.
The following statement (whose verification is now totally straightforward
and is left to the reader) makes precise the claim that this decoding
procedure correctly decodes on a tail of $U(\k)$:
\medskip

{\nar\smallskip\noindent $(\ast )$\ \ if $p \in \Cal I,\
\k \in d^p,\ \be^p(\k) \leq
\a < \de^p(\k)$, and $\a \in U(\k)$, then $G_1(\a) = g^p(\a)$.
\bigskip}

\noindent
{\bf (4.6.6)}
\medskip
Before turning to the remainder of the definition in the generic case, we
turn to the decoding of $g|[\k,\ \de^*)$
on the multiples of $\k$ only, from a $g$
defined on a large enough domain below $\k$, since in (6.1) (7) (b) we only
need this for the multiples of $\k$.  Here, rather than having $G$ defined on
$\bigcup \{(\la,\ \la^+)|\aleph_2 \leq \la = card\ \la < \k\}$, we have a tail
$t$ of $D_\k$ and a function $\de$ with domain $t$ such that for $\la \in t,\
\de(\la) \in U(\la)$ with $\la < \de(\la) < \la^+$, such that
$g$ is defined on $\bigcup \{(\la,\ \de(\la))|\la \in t\}$ (of course, $g$
will also be defined elsewhere, but only this is relevant for our decoding).
We take $\de^* = scale(\de)$.  Finally, to complete the description
of the context of (6.1) (7) (b), below, we have that there is very tidy
$p \in P$ such that $t = d^p \cap D_\k$ and that for $\la \in t,\
\de^p(\la) \leq \de(\la),\ g^p|(\la,\ \de^p(\la)) = g|(\la,\ \de^p(\la))$ and
that for $\a \in U(\la) \cap [\de^p(\la),\ \de(\la))$ we will have
$g|s_\a$ coding the $<_{L[A \cap \k]}\text{-least}$ well-ordering of $\la$ in
type $\a$ and $g(\a) =\ ?$ except possibly in the case of $\a = \de^p(\la)$
when it is also possible that $g(\de^p(\la)) =\ !$.

Having described the context, the procedure is essentially identical to the
above, with the obvious notational analogies, so we limit ourselves to
describing the more substantial differences.  These deal only with the
definitions of the functions
analogous to the $H^\a$ and $\pi^\a_1$, above.
In both cases, we impose
that $h^\a(\la),\ \pi^\a_1(\la) \leq \de(\la)$ by taking them as defined to be
the min of $\de(\la)$ and the value defined as above.  This completes the
treatment of the singular non-generic case.
\bigskip
\noindent
{\bf (4.6.7)}
\medskip
We complete the description of the decoding procedure by returning to the
final phase of the singular generic case:  defining $G$ on the non-multiples
of $\k$ in $(\k,\ \k^+)$, detecting contamination, and revising the definition
of $G_1$ on the bounded initial segment of multiples of $\k$ where there is
contamination.  In several places, our argument will appeal to density
arguments from (6.1).  We should emphasize {\it that there is no
circularity here,} since we have already completed the portion of the argument
(the singular non-generic case) needed in (6.1) (7) (b).

So, let $\k < \eta < \k^+$ with $\eta$ not a multiple of $\k$.
We argue that there is a tail, $T$, of $X_\eta$
and an $i \in \{0,\ 1\}$ such that $1 - i \not\in G_1[T]$ and for
a cofinal set of $\a \in T,\ G_1(\a) = i$.  We first show that
$X_\eta \cap G_1^{-1}[\{0\}],\ X_\eta \cap G_1^{-1}[\{1\}]$ cannot both be
cofinal.  To this end, let $p \in \Cal I$ such that $\k \in d^p$ and
$\eta < \de^p(\k)$ (such a $p$ exists, by (6.1) (4) and (7)).  Suppose, now,
towards a contradiction, that $\be^p(\k) \leq \a_0,\ \a_1 \in X_\eta$ and
that $G_1(\a_j) = j$.  Let $i\ :=\ g^p(\eta)$.  By (6.1) (4), again,
there are $q,\ r \in \Cal I$ with $\de^q(\k) > \a_0,\ \de^r(\k) > \a_1$.
Clearly we can assume that $p \leq q \leq r$.  By $(\ast )$, above,
$g^r(\a_j) = j$.  But this contradicts (3.4) (A), since $g^r(\a_j)$ is forced
to be $i$.  Finally, from $(\ast )$, above and (6.1) (4), it is immediate that
if $\a < \k^+$ there is $\a < \a^\prime$ and
$p \leq q \in \Cal I$ with $\a^\prime \in X_\eta \cap \de^q(\k)$ such that
$g^q(\a^\prime) = i$ and $\be^q(\k) \leq \a^\prime$.  Then, by
$(\ast )$, above,
again, $G_1(\a^\prime) = i$ and we are finished.

So, we define $G(\eta)\ :=$ that $i \in \{0,\ 1\}$ such that for a cofinal
set of $\a \in X_\eta,\ G_1(\a) = i$.  Finally, for $\a \in U(\k)$,
we say that $\a$ is $G_1\text{-contaminated}$ iff
(following (3.3)) $G_1(\a) \neq\ !$ and
either there are cofinal subsets,
$Y_i \subseteq b_\a,\ i = 1,\ 2$, and $x_1 \neq x_2$ such that for $\xi \in
Y_i,\ G(\xi) = x_i$, or, if $\a$
is an odd multiple of $\k^2$, letting $\ga$ be such that
$\a \in X_\ga,\ G_1(\a)
\in \{0,\ 1\}$ but $G_1(\a) \neq G(\ga)$.
We let $\be\ :=\
sup\ \{\a|\a\text{ is } G_1\text{-contaminated}\}$.  It is totally
straightforward (and left to the reader to verify) that if $p \in \Cal I$ and
$\k \in d^p,\ \be = \de^p(\k)$.  Then we define $G(\a)$ for $\k < \a < \k^+,\
\a$ a multiple of $\k$, by setting $G(\a)\ :=\ G_1(\a)$ if $\be \leq \a$,
while for $\a \in U(\k) \cap \be$, we set
$G(\a)\ :=\ ?$ and we define $G|s_\a$ to code the
$<_{L[A \cap \k]}\text{-least}$ well-ordering of $\k$ in type $\a$.  This
completes the decoding procedure.

%

\subheading{\S 5.  STRATEGIC CLOSURE AND DISTRIBUTIVITY}
\bigskip
\noindent
{\bf (5.1)}\ \ THE WINNING STRATEGY.
\medskip
In this item we prove
Lemma 5 of the Introduction.  So, let $\th \geq \aleph_2$ be
regular and fix $p_0 \in P_\th$, $\Cal M$ and $\orh{\Cal N}$ as in (1.1)
Further assume that $\orh{\Cal N}$ is super
$\Cal M$ coherent, and recall the definition of the game
$G(\th,\orh{\Cal N},p_0)$ in (1.1) (where $\bold Q = \bold P$ and
$X$ is the class of very tidy conditions).
Recall that BAD must play very tidy conditions.

GOOD's strategy will be to use (4.3) at successor stages, so that, in
the notation of (4.3), she will have
$p_{2\a + 2}\ :=\ q(p_{2\a + 1},\N)$, for an $\N$ which we shall
describe below.  This will be chosen so as to guarantee that
at limit stages, we have the hypotheses of (4.5), so that,
at limit stages, $\si$, GOOD will take $p_\si$ to be given by (4.5).
The other implicit assumption is that
at all stages so far, BAD has succeeded in
\lq\lq catching" $p_{2i},\ p_{2i + 1}$
inside $|\Cal N_{\a(i)}|$.

For the successor step, we take $\N^\prime\ :=\ \N_{\a(i)},\
p\ :=\ p_{2\a + 1},\ \N\ :=\ $ the Skolem hull in
$\M$ of $N^\prime \cup \{\N^\prime\}$
and $q\ :=\ q(p,\N)$.  Note that we easily
have that $p\chi_{\N^\prime} \in N$, so that for all s-like $\k \in d(p),\
\de^q(\k) > p\chi_{\N^\prime} > \de^p(\k)$.  This guarantees that
we have $g^q(p\chi_{\N^\prime}(\k))\ =\ !$.

Now, let $\si \leq \theta$ be a limit ordinal and let
$\k \in d(p)$ be singular,
where $p$ is as in (4.5).  Let $i_0$ be the least $i < \si$ such that
$\k \in d^{p_i}$.  Since $p_{i_0} \in |N_{\a(j_0)}|$, where
$j_0$ is least such that $i \leq 2j + 1,\ \k$ is $\Cal N_{\a(j)}\text{-
controlled}$ for all $j_0 \leq j < \si$.
Then, letting $\de\ :=\ \de^p(\k),\
\{p\chi_{\N_{\a(j)}}(\k)|j_0 \leq j < \si\}$
is a subset of $(g^p)^{-1}[\{!\}]$ which is cofinal in
$\de$.  Finally, the supercoherence of the model sequence $\orh{\Cal N}$
guarantees that it is also a subset of $C_\de$, which
means that it is a subset of $Z_\de$.  It is then routine to see that
we have the hypotheses of (4.5), and therefore, that the strategy
for GOOD is winning.  To obtain the distributivity properties,
we must see that BAD needn't lose due to inability to
\lq\lq catch" $p_{2i},\ p_{2i + 1}$
inside $|\Cal N_{\a(i)}|$, and, more importantly, that there are enough
supercoherent sequences.  These points are addressed in the next item;
the latter draws on the work of our companion paper \cite{17}.
\bigskip
\noindent
\proclaim{(5.2)  Corollary}
$\bold P_\th$ is $(\th,\ \infty)$-distributive.
\endproclaim
\demo{Proof}  If $p_0 \in P_\th$ and
$(\tilde D_i|i < \th)$ is a definable-in-parameters
sequence of open dense subclasses of $P_\th$, begin by picking singular $\nu$
with $cf\ \nu >> \th$, such that all parameters in the definition of
$(\tilde D_i|i < \th)$ lie in $H_\nu$ and such that $(H_\nu,\ \in)$ reflects
all $\Sigma_n$-formulas, where the definitions of $(\tilde D_i|i < \th)$ and
$\bold P_\th$ are $\Sigma_n$ and $n$ is larger than the number of quarks
in the physical universe.  Let $D_i = \tilde D_i \cap H_\nu$.  We take $\M =
(H_{\nu^+},\ \in,\ \cdots )$, we let $\N_0 \prec \M$, with
$\th + 1 \cup \{(D_i|i < \th), \ \bold P_\th|H_\nu \} \subseteq
N_0\ (:=\ |\N_0|)$ and $card\ N_0 = \th,\ [N_0]^{<\ \th} \subseteq N_0$.
By the main result of our companion paper,
\cite{17}, see (1.4),
we can find super $\M$-coherent
$(\N_i|i \leq \th)$ starting from $\N_0$.
We then play a run of the game $G(\th,\ \orh{\N},\ p_0)$,
where GOOD plays by her winning strategy.

We argue that BAD can produce a subsequence of $\orh{\N}$ which
\lq\lq catches" $p_{2i},\ p_{2i + 1}$ inside $N_{\a(i)}$.
For successor $i$, given that $p_{2i - 2},\ p_{2i - 1} \in N_{\a(i - 1)}$,
the last sentence of Remark 1 immediately gives that if BAD chooses
$\a(i) > \a(i - 1),\ p_{2i} \in N_{\a(i)}$.  Then, if he chooses
$p_{2i + 1}$ from $N_{\a(i)}$, this is as required.
For limit $i$, letting $\a^*\ :=\ sup\ \{\a(j)|j < i\}$, given that
all the $p_j \in N_{\a^*},\ j < i$, the Remark of (4.5) immediately
gives that if BAD chooses $\a(i) > \a^*$, then $p_i$, the $p$ of
(4.5), will lie in $N_{\a(i)}$.  Once again, if he then chooses
$p_{i + 1}$ from $N_{\a(i)}$, this will be required.

Thus, in a such a play,
we will actually produce a $p_\theta$.  If, in addition to the above,
BAD chooses $p_{2i + 1} \in D_i \cap N_{\a(i)}$, then clearly we will
have $p_\theta \in \bigcap \{D_i|i < \theta\}$, as required.
\enddemo
\medskip
\noindent
\proclaim{(5.3)  Remarks}
\roster
\item"{(1)}"  It follows from (5.1) and (5.2) that any iteration of
${\bold P}_\th$ with supports of cardinality $\leq \th$
is $(\th,\infty)$-distributive.  The reason is quite simply that
(5.1), (5.2) and the results of \cite{17}  depend only on
ground model properties.  In fact, S. Friedman has informed us that
prior to our work similar observations had been made about Jensen's
original coding and the variants of it mentioned in the Introduction.
\medskip
\item"{(2)}"  As mentioned in the Introduction, by a simple variant of the
proof of (5.1), we can show that GOOD has a winning strategy for the version
of the game where we start from a $p_0 \in P$, which is not necessarily $\in
P_\th$.  We \lq\lq freeze'' $g^{p_0}|\th$, (and therefore, also
$\be^{p_0}|\th$) and never add any $(\nu,\ \xi)$, with
$\xi < \th$, to any $\Xi(p_i)$.  We use this in (6.1) (7) (c).
\endroster
\endproclaim

%

\subheading{\S6. EXTENDABILITY AND CHAIN CONDITION}
\noindent
\proclaim{(6.1)\ \ Lemma}
\medskip
Let $\k > \aleph_1$ be a cardinal, $p \in P$.  Then, in
each of the cases $1 \leq i \leq 7$, below, there is a $q \in P,\ p \leq q$,
satisfying the conclusion of (i) (which follows the colon, in each case).
\endproclaim
\roster
\item $\k$ is a successor cardinal, $\k \not\in d^p: \k \in d^q$,
\item $\k$ is a regular cardinal, $\k \in d^p,\ \k \leq \a < \de^p(\k),
\ g^p(\a) = 1,\ \a \not\in W(p)$,
and if $\k$ is inaccessible, then letting
$\be_1$ be as in (3.4) (D), either $\be^p(\k) \leq \a <
\de^p(\k)$ or $\be_1 \leq \a < \k^2:  \a \in W(q)$,
\item $\k$ is regular, $\k \in d^p,
\ \de^p(\k) \leq \a < \k^+:  \a <
\de^q(\k)$,
\item $\k$ is singular, $\k \in d^p,\ \k < \ga < \de^p(\k),\ \ga$ is not a
multiple of $\k,\ \de^p(\k) < \a < \k^+$:  there is $\a^\prime \in X_\ga,
\ \a < \a^\prime < \de^q(\k)$ and $g^q(\a^\prime) = g^p(\ga)$,
\item $\k$ is regular, $\k \in d^p,\ \k \leq \a < \de^p(\k),
\ g^p(\a)
= 0,\ \sigma < \k$ and if $\k = \la^+$, where $\la$ is a limit cardinal,
then $\la \in
d^p$: there is $\zeta \in b_\a \setminus \sigma$ with $card\ \zeta \in
d^p,\ \be^q(card\ \zeta) \leq \zeta < \de^q(card\ \zeta)$ and
$g^q(\zeta) = 0$,
\item $\k$ is inaccessible, $\k \not\in d^p: \k \in d^q$,
\item $\k$ is singular, $\k \not\in d^p: \k \in d^q$.
\endroster

\demo{Proof} The properties are given in order of increasing difficulty
and/or dependence on earlier properties.  The analogue of (5), where $\k =
\la^+,\ \la$ a limit cardinal, $\la \not\in d^p$,
is achieved by first adding $\la$ to
$d^p$, via (6) or (7), as appropriate, then using (5).
We deal with the cases in the given order.
In all cases, in virtue of (4.3), we can assume that $p$ is very tidy.
In all cases except (5),
when we have to define $g^q(\a)$, we shall always make $g^q(\a) = 1$
unless it is forced to be $0$, in which case we make it $0$, {\bf except}
in the following cases:
\medskip

\roster
\item"{(a)}"\ \ $card\ \a$ is s-like and $\a$ is a multiple of $card\ \a$,
\item"{(b)}"\ \ $\k = card\ \a$ is inaccessible, $\a < \be^q(\k)$ and
$\a$ is a multiple of $\k$.
\endroster
\medskip
\noindent
Thus, in what follows, except in case (5)
we shall limit ourselves to treating the above cases.
In case (5), we will find a $\zeta$ {\bf which is not forced to be 1 and
which is not in} $\bold {R(p)}$ and for this $\zeta$, we shall make
$g^p(\zeta) = 0$.  $\zeta$ will not be a multiple of $card\ \zeta$.  This
will be the only exception to our general procedure, even in case (5).

For (1), if $\k = \la^+$, we simply set $d^q =d^p \cup \{\k\},\ \vX (q) =
\vX (p)$, we set $\de^q(\k) = \k^2,
\ \be^q(\k) = \k + 1$.

For (2), we shall have $g^q = g^p$.
If $\k = \la^+$, where $\la \not\in d^p$, let
$\xi = \la$.  If $\k = \la^+$, where $\la \in d^p$, let $\xi =
\de^p(\la)$.  If $\k$ is inaccessible and $\k \cap d^p = \emptyset$,
let $\xi = \aleph_2$.  Finally, if $\k$ is inaccessible and
$\k \cap d^p \neq \emptyset$, let $\xi = sup(\k \cap dom\ g^p)$.  In the last
two cases, if $\theta < \k,\ \theta$ is regular, we can
always take $\xi \geq \theta$, as well.  Then, we let $\vX (q) =
\vX (p) \cup \{(\a,\ \xi)\}$.

For (3), we will have $d^q = d^p,\ \be^q = \be^p,\ \vX (q) = \vX (p)$,
and for all $\mu \in d^p \setminus \{\k\},\ \de^q(\mu) = \de^p(\mu)$.
We set $\de^q(\k) =\ $ the least $\theta > \a$ which is a multiple of $\k^2$.
Now,
suppose $\de^p(\k) \leq \xi < \de^q(\k)$.  The only case we have to
treat is when $\k$ is s-like and regular, and
$\xi$ is a multiple of $\k^2$; we set $g^q(\xi) =\ ?$, and we take
$g^q|s_\xi$ to be as required by (3.4)(E).  Clearly $q$ is as required.

For (4), pick $\a^\prime \in
X_\ga,\ \a^\prime > \a$.  We shall have $d^q = d^p,\ \be^q = \be^p,
\ \vX (q) = \vX (p)$.  If $\mu \in d^p,\ \mu \not\in \{\k\} \cup \Delta_\k$,
we shall also have $\de^q(\mu) = \de^p(\mu)$.

We set $\de^q(\k) = \a^\prime + \k^2$.
For a tail of $\la \in D_\k,\ \la \in
d^p$ and $\de^p(\la) \leq f^*_{\de^p(\k)}(\la) < f^*_{\a^\prime}(\la) <
f^*_{\de^q(\k)}(\la)$.  So let $\la_0 = \la_0(\a^\prime)$ be sufficiently
large so that the preceding holds for $\la \in D_\k \setminus \la_0$.
Clearly, we may assume $\la_0 \not\in D_\k$.  For
$\la \in D_\k \setminus \la_0$,
let $\de^q(\la) = f^*_{\de^q(\k)}(\la)$.  If $\k$ is a limit of singular
cardinals and $\tau > \la_0$ is a successor point of $D_\k$, let $\nu\ :=\
\de^q(\tau)$ and let $\nu^\prime = f^*_{\a^\prime}(\tau)$.
Then, as for $\k$ and $\de^p(\k)$, there is
$\la_0(\nu^\prime) \geq \la_0$ such that for all $\la_0 \leq \la \in D_\tau,\
\de^p(\la) \leq f^*_{\de^p(\tau)}(\la) < f^*_{\nu^\prime}(\la) <
f^*_\nu(\la)$.  For such $\la$ we set $\de^q(\la)\ :=\ f^*_\nu(\la)$.
For all other $\la \in \Delta_\k
\cap d^p,\ \de^q(\la) = \de^p(\la)$.

Suppose now that $\la \in d^p$ and $\de^p(\la) < \de^q(\la)$ (so, in
particular $\la \in \{\k\} \cup \Delta_\k$).
We deal first with defining $g^q(\de^p(\la))$, so let $\xi = \de^p(\la)$.
If $\xi$ is $p\text{-deactivated, }$ we set $g^q(\xi)\ =\ !$;
otherwise, we set $g^q(\xi)\ =\ ?$.  In both cases we take
$h^q(\xi) = \xi$, and define $g^q|s_\xi$ to satisfy (3.4) (E).
In virtue of the rest of the definition of $g^q$, Remark 2 of (4.3)
will guarantee that this is as required (recall that $p$ is very tidy!).

Next, suppose that $\xi$ is a multiple of $\la^2$ with $\de^p(\la) < \xi
< \de^q(\la)$.  We shall have that $g^q(\xi)\ =\ ?$ {\bf unless} one of
the following occurs:

\medskip
\roster
\item"{(c)}"\ \ $\la = \k$ and $\xi = \a^\prime$,

\item"{(d)}"\ \ $\la \in D_\k \setminus \la_0(\a^\prime)$ and
$\xi = f^*_{\a^\prime}(\la)$,

\item"{(e)}"\ \ $\k$ is a limit of singular cardinals and there is $\tau$,
a successor point of $D_\k \setminus \la_0(\a^\prime)$ such that,
letting $\nu^\prime :=\ f^*_{\a^\prime}(\tau),\ \la \in D_\tau \setminus
\la_0(\nu^\prime)$ and $\xi = f^*_{\nu^\prime}(\la)$.
\endroster
\medskip

\noindent
In these cases, we set $g^q(\xi) = g^p(\ga)$.
In all cases, we will have $h^q(\xi) = \xi$, and we'll define
$g^q|s_\xi$ to satisfy (3.4) (E).
By Remark 2 of (4.3), and the fact that if $\la \in d^p$ is singular
and $\de^p(\la) < \nu$ then $f^*_\nu >^* \de^p|D_\la$, it is then clear
that $q \in P$ and is as required.  In fact, it is easily verified that
$q$ is very tidy, though we do not need this.

For (5), suppose, first, that $\k = \la^+$, and $\la \in d^p$.  In this
case, we shall have
$d^q = d^p,\ \be^q = \be^p,\ \vX (q) = \vX (p)$.  By (3.2)(A)
for $p$ and $\k^\prime = \k$, we can find $\zeta_0 < \k$ such that
whenever $(\eta,\ \xi) \in \Xi(p)$ and $\xi < \k$,
then $b_\a \cap b_\eta \subseteq
\zeta_0$.  Without loss of generality, $\zeta_0 < \de^p(\lambda)$.
Now let $\zeta \in
b_\a \setminus max\ (\sigma,\ \zeta_0)$.
By our choice of $\zeta_0,\ \zeta \not\in R(p)$.  Also,
since $\zeta \in b_\a$ and $\a \in (\k,\ \k^+)$ where $\k$ is a successor
cardinal, $\zeta$ is even, but not a multiple of $\la$.  Thus,
$\zeta$ is not forced to be 1.  Accordingly, we set $g^q(\zeta) = 0$.
The remainder of the construction of $q$ divides into
cases, according to whether $\la$ is regular or singular.

If $\la$ is regular, we
proceed as in (3), with $\la$ in place of $\k$ and $\zeta$ in place
of $\a$  and with the already-noted difference that $g^q(\zeta) = 0$.
If $\la$ is singular, we proceed as in (4), with $\la$
in place of $\k,\ \zeta$ in place of $\a$ with the already-noted
difference that $g^q(\zeta) = 0$;  the argument here is simpler than
in (4) since there are no $\ga$ nor $\a^\prime$ involved.
Then $q$ is as required.

If $\k = \la^+,\ \la \not\in d^p$, then, by hypothesis, $\la$ is a successor
cardinal, so we can use case (1) to obtain $p \leq q^\prime$ with $\la \in
d^{q^\prime}$, and then apply the immediately preceding argument to $q^\prime$
instead of $p$, to obtain the required $q$.  Finally, suppose $\k$ is
inaccessible.  As above, we can find $\zeta_0 < \k$ such that whenever
$(\eta,\ \xi) \in \Xi(p)$
and $\xi < \k$, then
$b_a \cap b_\eta \subseteq \zeta_0$.  Pick
$\k^\prime \geq max(\sigma,\ \zeta_0),\ \k^\prime = \al_\tau,\ \tau$
even successor.  We take $\zeta\ :=\ f^*_\a(\k^\prime)$.
We note, once again, that by our choice of $\zeta_0,\ \zeta \not\in R(p)$, and
that since $\zeta$ is even and $card\ \zeta$ is an even successor,
$\zeta$ is not forced to be 1.
Then, we can proceed as in (1) and (3), to add $\k^\prime$ to $d^p$,
and make $\de^q(\k^\prime) > \zeta$,
\underbar{EXCEPT} that, as above, we can also make
$g^q(\zeta) = 0$.  Clearly $q$ is as required.

For (6), we shall have $d^q = d^p \cup \{\k\},\ \vX (q) = \vX (p)$ and
for $\mu \in d^p,\ \de^q(\mu) = \de^p(\mu)$.
By (3.2)(A), for p with $\k^\prime =
\k,\ d = d^p$, we can compute $\be^q(\k)$ according to (3.4)(C) and
$\be_1$, according to (3.4)(D) and they will be bounded,
in $\k^+,\ \k^2$, respectively.
We take $\de^q(\k) =\ $ the
least multiple of $\k \geq \be^q(\k)$.
If $\k < \a < \be_1$ or
$\k^2 \leq \a < \be^q(\k)$ and $\a$ is a multiple of $\k$, we set
$g^q(\a) =\ ?$ iff $\a$ is contaminated.
If it is not contaminated, we set $g^q(\a) = 1$, unless it is forced
to be 0, in which case we make $g^q(\a) = 0$.
If $\be^q(\k) \leq \a < \de^q(\k),\ \a$ is a multiple of $\k$, we set
$g^q(\a) = 1$ unless it is forced to be
$0$; in this case, we set $g^q(\a) = 0$.  Finally, we define
$g^q$ on the multiples of $\k$ in $[\be_1,\ \k^2)$ to
satisfy (3.4) (D), (1) and (2).  Clearly this $q$ is as required.

Case (7) divides into subcases, as follows:

(a) $D_\k \cap d^p$ is bounded in $\k$ (in some sense, the simplest subcase:
we must add to $d^p$ a tail of $\Delta_\k$, but there is no contamination),

(b)  $D_\k \subseteq^* d^p$ ($\k$ will be the only new member of $d^q$,
but we must deal with contamination); we shall use the decoding procedure of
(4.6);

(c)  (a) and (b) both fail (the most complicated case:  we must combine the
methods used for (a) and (b), \underbar{and} appeal to (5.1), (5.2)).

In case (a),
let $\la_0 < \k$ be such that
$D_\k \cap d^p \subseteq \la_0$.  Clearly, we
may assume $\la_0 \not\in D_\k$, and, anticipating the argument for (c), if
$\theta < \k,\ \theta$ regular, we can take $\la_0 \geq \theta$.  We shall have
$\vX (q) = \vX (p),\ d^q = d^p \cup \{\k\} \cup
(\Delta_\k \setminus \la_0)$.
For $\mu \in d^p$,
we will have $\de^q(\mu) = \de^p(\mu)$.  For
$\la \in \{\k\} \cup (\Delta_\k \setminus \la_0)$, we shall
have $\be^q(\la) = \la + 1$ if $\la$ is a successor cardinal and
$\be^q(\la) = \la^2$, if $\la$ is singular.

We set $\de^q(\k) = \k^2$.  For $\la \in D_\k \setminus \la_0$,
we set $\de^q(\la) = f^*_{\k^2}(\la)$.  If $\k$ is a limit
of singular cardinals, $\la$ a successor point of $D_\k,\ \la > \la_0$, let
$\de = \de^q(\la)$.  Then, if $\tau \in D_\la \setminus \la_0$, we set
$\de^q(\tau) = f^*_\de(\tau)$.  Then, for $\la \in \{\k\} \cup
(\Delta_\k \setminus \la_0)$, if $\la < \a < \de^q(\la)$, and $\a$ is
a multiple of
$\k^2$, we set $g^q(\a) =\ ?$, we take
$h^q(\a) =
\a$, as required by (3.4) (E), and we define
$g^q|s_\a$ to code this, as required by (3.4) (E) (we can always find
$R \in L[A \cap \la]$ as required, since either $\la$ is singular, in which
case $(\la^+)^L = \la^+$, or $\la \not\in \Lambda$, in which case
$(\la^+)^{L[A \cap \la]} = \la^+$).  If $\la$ is s-like and
regular, recall that $f^*_\a$ was defined at the end of (1.2).
This completes the proof in case (a) of (7).

In case (b), we will have $d^q = d^p \cup \{\k\},\ \vX (q) = \vX (p)$ and
for $\mu \in d^p,\ \be^q(\mu) = \be^p(\mu)$.  If $\mu \in d^p,\ \mu
\not\in \Delta_\k$, we shall also have $\de^q(\mu) = \de^p(\mu)$.

We set $\de^q(\k) = \de^* = scale(\de^p|D_\k)$.
Let $\la_0 < \k$ be such that if $\la \in
D_\k \setminus \la_0$, then $\la \in d^p\ \&
\ \de^p(\la) \leq f^*_{\de^*}(\la)$.  Clearly, we
may assume $\la_0 \not\in \Delta_\k$ and, anticipating the
argument for (c) when $\k$ is a limit of singular cardinals, below, if
$\theta$ is regular, $\th < \k$ we can take $\la_0 \geq \th$.  For $\la \in
D_\k \setminus \la_0$ we set $\de^q(\la) = f^*_{\de^*}(\la)$.
If $\k$ is a limit of singular cardinals, for
such $\la$, if $\la$ is a successor point of $D_\k$ and $\de_0 = \de^p(\la) <
\de^q(\la) = \de_1$, then, on a tail of $\eta \in D_\la,\ \eta \in d^p$ and
$\de^p(\eta) = f^*_{\de_0}(\eta) < f^*_{\de_1}(\eta)$.  So, let $\eta_0 =
\eta_0(\la)$ be such that whenever
$\eta_0 \leq \eta \in D_\la,\ \eta \in d^p$
and $\de^p(\eta) = f^*_{\de_0}(\eta) < f^*_{\de_1}(\eta)$.  For such $\eta$,
set $\de^q(\eta) = f^*_{\de_1}(\eta)$.
If $\tau \in D_\k \cap \la_0$, or (if $\k$ is a limit of singular
cardinals), for some successor point, $\la$, of
$D_\k \cap \la_0,\ \tau \in D_\la$, or (if $\k$ is
a limit of singular cardinals) for some
successor point, $\la$, of $D_\k \setminus \la_0,\ \tau
\in D_\la \cap \eta_0(\la)$, set $\de^q(\tau) = \de^p(\tau)$.

For $\la \in \Delta_\k$ such that $\de^p(\la) < \de^q(\la)$,
we handle the definition of
$g^q|[\de^p(\la),\ \de^q(\la))$ as we did in case (4), except that,
here again, as in (5), the argument is simpler since there are no
$\ga,\ \a^\prime$ involved.

Thus, it remains to define $\be^q(\k)$ and $g^q|(\k,\de^q(\k))$.  We define
$g^q$ in the usual way on the non-multiples of $\k$ in $(\k,\ \de^*)$.
We shall define $\be^q(\k)$ to satisfy (3.4)(C)
with $d = d^p$ and $q|\k$ in
place of $p|\k$.
Note that conceivably $\overline\de < \beta^q(\k) < \de^*$, where
$\overline\de$ is the least multiple of $\k^2,\ \de$,
with $\k < \de \leq \de^*$ such that
$\neg(f^*_\de \leq^* \de^p|D_\k)$, since instances of contamination
could arise due to the definition of the $g^q(\la)$ for those $\la \in
D_\k \setminus \la_0$ with $\de^p(\la) < \de^q(\la)$ (if there
are cofinally many such).

For $\a$ a multiple of $\k^2,\ \k < \a < \be^q(\k)$, we set
$g^q(\a) =\ ?$, and we define
$g^q|s_\a$ to satisfy (3.4) (E).  For $\be^q(\k) \leq \a < \de^*,\ \a$ a
multiple of $\k^2$, we define $g^q|(\{\a\} \cup s_\a)$ by recursion on
$\a$, following the singular non-generic case of the decoding procedure
of (4.6), with $g = g^q|\k$.  This completes
the construction of $q$ in case (b).

For case (c), our strategy is to obtain $p \leq p^\prime,
\ p^\prime \in P$ such that
the hypothesis of case (a), above, holds for $p^\prime$ and $\k$, and then
apply (a) to $p^\prime$.
The construction of $p^\prime$ differs according to whether
$\k = \la^{+\omega}$
or $\k$ is a limit of singular cardinals.  The former case is \underbar{much}
easier, and we consider it first.  Here, we obtain $p^\prime$
by \underbar{simultaneously} adding $\la$ to $d^p$, following the procedure
of (1), for $\la \in$ any final segment of $D_\k \setminus d^p$.  In
particular, anticipating the argument when $\k$ is a limit of
singular cardinals, the final segment can be taken to lie above $\theta$, if
$\theta$ is regular, $\theta < \k$.  The simultaneity is emphasized to make
clear that we \underbar{are not yet} appealing to any
strategic closure properties.  We then proceed as in (a) with
$p^\prime$ in place of $p$.

When $\k$ is a limit of singular cardinals, as a first step toward
obtaining the desired $p^\prime$, we first
simultaneously add to $d^p$ all the $\tau \in D_\la$, for $\la$ a successor
point of $D_\k,\ \tau \not\in d^p$,
according to the procedure for (1).  This
is a condition, $p_0$, intermediate between $p$ and $p^\prime$.

To obtain $p^\prime$, we let $(\la_i: i < \sigma)$ increasingly enumerate
$D_\k \setminus d^p$.  We let $\theta > \aleph_2$ be regular,
$\sigma \leq \theta < \k$.

Let $\Cal M$ be a master model, $\Cal M = (H_{\nu^+},\ \in,\ \cdots),
\ \nu$ singular, $\nu >> \k$, such that $p_0 \in H_\nu$ and
$(H_\nu,\ \in)$ models a sufficiently rich fragment of ZFC, etc., as in
(5.1).  As in (5.1), we can assume that we have $(\Cal N_i: i \leq \theta)$
which is super $\Cal M\text{-coherent}$, with $p_0 \in |\Cal N_0|$.  So,
fix such $(\Cal N_i: i \leq \theta)$.
We then play the
following run of the variant of $G(\theta,\ \orh{\Cal N},\ p_0)$, mentioned in
(5.3).  GOOD plays by the
winning strategy of (5.1).
BAD chooses $\a(i)$ as in (5.2), and
obtains $p_{2i + 1} \in
|\Cal N_{\a(i)}|$, by adding
$\la_i$ to $d^{p_{2i}}$, following the procedure of case (a) of (7)
(note that the hypothesis of
case (a) will always hold for $\la_i$ and $p_{2i})$.  Then,
$p^\prime$ can be taken to be $p_{\sigma}$.  This completes the proof for (c),
when $\k$ is a limit of singular cardinals, and therefore completes the proof
of (7) and the Lemma.
\enddemo
\bigskip
\noindent
{\bf (6.2)}\ \ $\bold {\dot {P}^\theta}$ HAS THE $\theta^+$-CHAIN CONDITION
\medskip
Let $\theta > \aleph_2$ be regular.  The crucial observation is:

\proclaim{(6.2.1) Proposition}  Suppose $p,\ q \in P,\ (p)_\theta,\ (q)_\theta$
are compatible in $\bold P_\theta,\ g^p|\theta =
g^q|\theta$ and $\be^p|\theta = \be^q|\theta$.
Then $p,\ q$ are compatible in $P$.
\endproclaim

\demo{Proof} Let $r \in P_\theta$ with $(p)_\theta,\ (q)_\theta \leq r$.  Note
that, without loss of generality, we may assume that $W(r) =
W_\theta(p) \cup W_\theta(q)$.
We shall show that $r^* \in P,\ p,\ q \leq r^*$, where
$r^* = (g^r \cup g^p|\theta,\ \be^r \cup \be^p|\theta,\ \vX (r^*))$,
where $\vX (r^*) = \vX(p) \cup \vX(q) \cup \vX(r)$.
Of course, $p,\ q \leq r^*$
is clear, once we've verified that $r^* \in P$.

For this, all clauses of (2.2), (3.2)
are clear, as are (3.4) (E) and all clauses of (3.5).  We argue
that there is no new contamination in $r^*$, from which it will follow readily
that we also have all (3.4) (A) - (D).  This will complete the proof.
Clearly there is no new contamination at singulars, and there is no new
contamination at inaccessibles above $\theta$.  So, suppose that $\k$ is
inaccessible, $\k \leq \theta$.  Suppose that $\a \in (\k,\ \k^+)$ and that
$\a$ is contaminated by $\a^\prime$.
Let $(\a^\prime,\ \xi) \in \vX (r^*)$ witness this, as in (3.3).
Then, $\a^\prime \in W(p) \cup W(q)$,
and since $\xi < \k \leq \theta$
clearly $(\a^\prime, \ \xi) \in \{\Xi(p),\ \Xi(q)\}$.
But then by the hypotheses of the Proposition, $\a$ must be contaminated by
$\a^\prime$ either in $p$ or in $q$ according to whether
$(\a^\prime,\ \xi) \in \Xi(p)$
or $\in\ \Xi(q)$.  This completes
the proof.
\enddemo

\proclaim{(6.2.2) Corollary}  In $V^{\bold {P_\theta}},
\ \bold {\dot {P}^\theta}$ has the
$\theta^+$-chain-condition.
\endproclaim

\demo{Proof} This is clear from (6.2.1) and the easy computation that
$\{(g^p|\theta,\ \be^p|\theta)|p \in P\}$ has power $\theta$, for all
regular $\theta > \aleph_2$.
\enddemo

\enddocument